\documentclass[10pt]{article}
\usepackage[all]{xy}
\usepackage{amsmath}
\usepackage{amsfonts}
\usepackage{amssymb}
\usepackage{amscd}
\usepackage{amsthm}
\usepackage{latexsym}
\usepackage{amsbsy}
\newcommand{\A}{\mathcal{A}}
\newtheorem{lemma}{Lemma}[section]
\newtheorem{proposition}{Proposition}[section]
\newtheorem{theorem}{Theorem}[section]
\newtheorem{corollary}{Corollary}[section]
\newtheorem{definition}{Definition}[section]
\newtheorem{example}{Example}[section]

\newtheorem{remark}{Remark}[section]

\newcommand{\nc}[2]{\newcommand{#1}{#2}}
\newcommand{\rnc}[2]{\renewcommand{#1}{#2}}
\rnc{\theequation}{\thesection.\arabic{equation}}
\nc{\beq}{\begin{equation}}
\nc{\eeq}{\end{equation}}
\rnc{\[}{\beq}
\rnc{\]}{\eeq}
\nc{\ba}{\begin{array}}
\nc{\ea}{\end{array}}
\nc{\bea}{\begin{eqnarray}}
\nc{\beas}{\begin{eqnarray*}}
\nc{\eeas}{\end{eqnarray*}}
\nc{\eea}{\end{eqnarray}}
\nc{\be}{\begin{enumerate}}
\nc{\ee}{\end{enumerate}}
\nc{\bd}{\begin{diagram}}
\nc{\ed}{\end{diagram}}
\nc{\bi}{\begin{itemize}}
\nc{\ei}{\end{itemize}}
\nc{\bpr}{\begin{proposition}}
\nc{\bth}{\begin{theorem}}
\nc{\ble}{\begin{lemma}}
\nc{\bco}{\begin{corollary}}
\nc{\bre}{\begin{remark}}
\nc{\bex}{\begin{example}}
\nc{\bde}{\begin{definition}}
\nc{\ede}{\end{definition}}
\nc{\epr}{\end{proposition}}
\nc{\ethe}{\end{theorem}}
\nc{\ele}{\end{lemma}}
\nc{\eco}{\end{corollary}}
\nc{\ere}{\end{remark}}
\nc{\eex}{\end{example}}
\nc{\bpf}{{ Proof. ~~}}
\nc{\epf}{\hfill\mbox{$\square$}\vspace*{3mm}}
\nc{\hsp}{\hspace*}
\nc{\vsp}{\vspace*}
\nc{\mc}{\mathcal{C}}
\nc{\ws}{\widetilde{S}}
\nc{\wdel}{\widetilde{\delta}}
\nc{\wsigma}{\widetilde{\sigma}}
\nc{\wtau}{\widetilde{\tau}}
\nc{\wvfy}{\widetilde{\varphi}}
\nc{\ra}{\rightarrow} \nc{\lra}{\longrightarrow}
\nc{\lla}{\longleftarrow}
\nc{\epn}{\varepsilon}
\nc{\si}{\psi}
\nc{\Del}{\Delta}
\nc{\del}{\delta}
\nc{\ro}{\rho}
\nc{\fy}{\phi}
\nc{\tu}{\tau}
\nc{\vfy}{\varphi}
\nc{\sig}{\sigma}
\nc{\alf}{\alpha}
\nc{\bt}{\beta}
\nc{\mh}{\mathcal{H}}
\nc{\ma}{\mathcal{A}}
\nc{\mbc}{\mathbb{C}}
\nc{\mbr}{\mathbb{R}}
\nc{\mbz}{\mathbb{Z}}
\nc{\mfg}{\mathfrak{g}}
\nc{\ug}{U(\mfg)}
\nc{\tg}{T(\mfg)}
\nc{\sg}{S(\mfg)}
\nc{\opl}{\oplus}
\nc{\bopl}{\bigoplus}
\nc{\Lam}{\Lambda}
\nc{\lam}{\lambda}
\nc{\weg}{\wedge}
\nc{\bweg}{\bigwedge}
\nc{\blt}{\bullet}
\nc{\ydh}{{}_{H}^{H}\mathcal{YD}}
\nc{\nn}{\nonumber}
\nc{\la}{\leftarrow}
\nc{\da}{\downarrow}
\nc{\ua}{\uparrow}
\nc{\ot}{\otimes}
\nc{\uh}{\underline{H}}
 \nc{\udel}{\underline{\Del}}
\nc{\us}{\underline{S}}
 \nc{\um}{\underline{m}}
\nc{\uepn}{\underline{\epn}}
\nc{\ueta}{\underline{\eta}}
\begin{document}

\title{Twisted Spectral Triples and Connes' Character Formula}
\author {\bf{Farzad Fathizadeh}\\
Department of Mathematics and Statistics\\
York University\\
Toronto, Ontario, Canada, M3J 1P3\\
ffathiza@mathstat.yorku.ca\\ 
\\
 \bf{Masoud Khalkhali}\\
 Department of Mathematics\\
 The University of Western Ontario\\
 London, Ontario, Canada, N6A 5B7\\
 masoud@uwo.ca
 }
\date{}
\maketitle
\begin{abstract}
We  give a proof of an 
analogue of Connes' Hochschild character 
theorem for twisted spectral triples obtained from twisting a spectral triple by scaling automorphisms, under some suitable conditions.  We also survey  some of the properties of twisted spectral triples that are known so far.

\end{abstract}

%\tableofcontents

\section{Introduction}
In \cite{conmos2} Connes and Moscovici introduced a refinement of the  notion of 
spectral triple called {\it  twisted spectral triple}. As is well known by now, the metric aspects
  of a geometric space in noncommutative geometry is encoded by
a {\it spectral triple} \cite{acbook}. This concept however puts  rather severe
restrictions on the underlying noncommutative space, as it forces the existence of a tracial 
state. In particular    it is 
inapplicable in \emph{type III situations} according to Murray-von Neumann's  classification of 
von Neumann algebras into types. It was to deal with these type III examples that twisted spectral triples were introduced. 

In a twisted spectral triple, the twist is afforded by an {\it automorphism}. More precisely, 
there is an automorphism $\sigma$ of the underlying noncommutative space $\mathcal{A}$
such that all twisted commutators $[D, a]_{\sigma}= D a- \sigma (a) D$ are bounded operators
($D$ is a selfadjoint operator playing the role of the Dirac operator).  While there are many natural examples of 
twisted spectral triples, proving general results about them, parallel to the untwisted 
case,  have proved to be a  challenging task. The most celebrated among such questions  are  perhaps
the extension of Connes' character formula for the Hochschild class of the Connes-Chern character of a spectral triple \cite{acbook},  
as well as an extension of Connes-Moscovici's local index formula, \cite{conmos1} to  twisted spectral
triples. In this paper we give a proof of a character formula for twisted spectral triples obtained from twisting a spectral triple by scaling 
automorphisms, under suitable conditions. We shall also briefly outline Moscovici's local index formula for twisted spectral triples obtained 
by conformal  perturbation of a spectral triple \cite{mos1}. 

This paper is organized as follows. In Section \ref{From spectral triples  to twisted spectral triples} 
we recall the basic notions of 
spectral triples, noncommutative integration and residue functionals for a spectral 
triple. We also recall Connes and Moscovici's local index formula \cite{conmos1}. We then recall
the notion of twisted spectral triples \cite{conmos2}. In Section \ref{Examples of twisted spectral triples} 
we recall several general methods to construct twisted spectral triples from 
\cite{conmos2}, \cite{mos1} and \cite{fatkha1}. An important idea here is twisting a spectral triple by 
scaling automorphisms and the corresponding algebra of
twisted pseudodifferential operators introduced by Moscovici \cite{mos1}. 
In Section \ref{Properties of twisted spectral triples} we recall some 
non-obvious properties of twisted
spectral triples from \cite{conmos2}, among them  the fact that one can define a Connes-Chern character 
for twisted spectral triples with 
values in ordinary non-twisted cyclic cohomology. Equivalently one can define a pairing between a 
finitely summable twisted spectral triple 
and $K$-theory. In section \ref{Hochschild class of the Connes-Chern character} we give a proof of an 
analogue of Connes' Hochschild character 
theorem for twisted spectral triples obtained from twisting a spectral triple by scaling automorphisms. Finally 
in the last section we recall Moscovici's Ansatz and proof of the local index formula for twisted 
spectral triples obtained by scaling automorphisms \cite{mos1}.

\section{From spectral triples to twisted spectral triples} \label{From spectral triples 
to twisted spectral triples}
The notion of a geometric space in noncommutative geometry is encoded by
a {\it spectral triple}. This concept however puts  rather severe
restrictions on the underlying noncommutative space that renders it
inapplicable in \emph{type III situations} according to Murray-von Neumann classification of 
von Neumann algebras into types. We shall first try to explain this point.

We start with a quick
review of the {\it Dixmier trace} and the {\it noncommutative integral},
following closely \cite{acbook}. Let $\mathcal{H}$ be a Hilbert space and let 
$\mathcal{K}(\mathcal{H})$ denote the two sided ideal
of compact operators on $\mathcal{H}$.
For a compact  operator $T: \mathcal{H} \to \mathcal{H}$,  let
$$ \mu_1(T) \geq \mu_2(T) \geq \cdots \geq 0$$
denote the sequence of eigenvalues
of its absolute value $|T|:=(T^*T)^{\frac{1}{2}}$,  written in  decreasing order. Thus, by the  
minimax principle, $\mu_1 (T) = ||T||$, and in general
$$ \mu_n (T) = \text{inf} \, \,||T|_V||, \quad n\geq1,$$
where the infimum is over the set of subspaces of codimension $n-1$, and
$T|_V$ denotes  the restriction of $T$ to the subspace $V$. A compact operator 
$T \in \mathcal{K}(\mathcal{H})$
is called trace class if $\sum \mu_n (T) < \infty$, and in that case the trace of $T$ is defined by
\[ \textnormal{Trace}(T) = \sum_{n} (Te_n, e_n) \nonumber \]
where $\{e_n\}_{n=1}^{\infty}$ is an orthonormal basis for $\mathcal{H}$. It is known that this
 is the unique
normal trace on the $C^*$-algebra $\mathcal{L}(\mathcal{H})$ of bounded operators on 
$\mathcal{H}$. The question of existence of non-normal traces on $\mathcal{L}(\mathcal{H})$ was left 
open until it was settled affirmatively by Jacques Dixmier. 
In \cite{dix} 
Dixmier shows that there are uncountably many non-normal traces on $\mathcal{L}(\mathcal{H})$. 
Many years later, Alain Connes discovered that these non-normal traces can in fact be used to define 
a process of noncommutative integration in noncommutative geometry as we describe next.

The natural domain of a Dixmier trace is the set of  operators
$$\mathcal{L}^{1,\infty } (\mathcal{H}) = \{T \in
\mathcal{K}(\mathcal{H}); \quad \sum_1^N \mu_n (T)=O \,(\text{log} N )\}.$$
Notice that trace class operators are automatically in $\mathcal{L}^{1,\infty} (\mathcal{H})$.
The Dixmier trace of an operator $T \in \mathcal{L}^{1,\infty}
(\mathcal{H})$
  measures the {\it logarithmic divergence} of
its ordinary trace. More precisely,  we are interested in the \emph{limit} of the bounded sequence
$$\sigma_N (T)= \frac{\sum_1^N \mu_n(T)}{\text{log} N}, \quad N=1, 2, \dots$$
as $N\to \infty$.
The first problem of course is that, while by our assumption the sequence is bounded,  the usual limit may
not exist and must be replaced by a \emph{generalized limit}, similar to \emph{Banach limits} of  non-convergent 
bounded sequences. A more challenging task  is to make sure that our  generalized limit still defines a trace.

To this end, let $\text{Trace}_{\Lambda} (T), \Lambda \in [1, \infty)$,  be
the piecewise affine interpolation of the partial trace function
$\text{Trace}_{N} (T) =\sum_1^N \mu_n (T)$.  Recall that a state on a  $C^*$-algebra
is a non-zero positive linear functional on the algebra. Let
$\omega : C_b[e, \infty) \to \mathbb{C}$ be a normalized state on the algebra of bounded 
continuous functions
on  $[e, \infty)$ such that $\omega (f)=0$ for all $f$ vanishing at $\infty$. Now, using $\omega$,  the
Dixmier trace of a positive operator $T \in \mathcal{L}^{1,\infty}
(\mathcal{H})$ is defined as
$$ \text{Tr}_{\omega}(T): = \omega (\tau_{\Lambda} (T)),$$
where
$$\tau_{\Lambda} (T)=  \frac{1}{\text{log} \, \Lambda} \int_e^{\Lambda} \frac{ 
\text{Trace}_{r} (T)}{\text{log}\,r}\frac{dr}{r}$$
 is the  Cezaro mean  of the function $\frac{ \text{Trace}_{r} (T)}{\text{log}r}$
 over the multiplicative group $\mathbb{R}^*_+$. One then extends
 $\text{Tr}_{\omega}$ to all of $\mathcal{L}^{1,\infty}
(\mathcal{H})$ by linearity.

The resulting linear functional  $\text{Tr}_{\omega}$ is a
  positive trace on $\mathcal{L}^{1, \infty} (\mathcal{H})$.  It is easy to see from its definition
  that if  $T$ actually
  happens to be a trace class operator then $\text{Tr}_{\omega}(T)=0$
  for all $\omega$, i.e., the Dixmier trace is invariant under
  perturbations by trace class operators. This is a very useful property and  makes $\text{Tr}_{\omega}$ 
  a  flexible tool in computations.
  
 The Dixmier trace,   $\text{Tr}_{\omega}$,  in general depends on the
limiting procedure $\omega$,   however,  for the class of operators
$T$  for which   $\text{Lim}_{\Lambda \to \infty}\,  \tau_{\Lambda}(T) $
exits, it  is  independent of the choice of $\omega$ and is equal to
$\text{Lim}_{\Lambda \to \infty} \tau_{\Lambda}(T)$. One of the  main results proved in
\cite{con3} is that if
$M$ is a closed  $n$-dimensional manifold,  $E$ is a smooth vector bundle on  $M$, $P$
is a  pseudodifferential
operator of order $-n$  acting between $L^2$-sections of E, and $\mathcal{H} = L^2 (M, E),$ then
$P\in \mathcal{L}^{1, \infty} (\mathcal{H})$   and, for any choice of
$\omega$,
$\text{Tr}_{\omega} (P)= n^{-1}
\text{Res} (P)$. Here Res denotes Wodzicki's  noncommutative residue \cite{wod1, kas1}.
For example, if $D$ is an elliptic first order differential
operator,  $|D|^{-n}$ is a pseudodifferential operator of order $-n$, and the Dixmier trace
$\text{Tr}_{\omega}(|D|^{-n})$ is independent of the choice of $\omega$.

Next, we would like to explain the notion of spectral triple. This concept has its roots in $K$-homology 
and Riemannian geometry simultaneously. We start by explaining the notion of a Fredholm module which 
is the conformal counterpart of a spectral triple.  

An \emph{odd Fredholm module} over a unital algebra $\mathcal{A}$ is a pair $(\mathcal{H}, F)$ 
where $\mathcal{H}$ is a Hilbert
space on which the algebra $\mathcal{A}$ acts by bounded operators and 
$F \in \mathcal{L}(\mathcal{H})$ is a selfadjoint
operator such that $F^2= id$, and such that the commutators $[F, \pi(a)]$ are compact 
operators for all $a \in \mathcal{A}$. Here $\pi: \mathcal{A}  \to  \mathcal{L}(\mathcal{H})$
denotes a unital action of $\mathcal{A}$ on $ \mathcal{H}$. A Fredholm module is called
 \emph{$p$-summable} $(1 \leq p < \infty)$,
if $[F, \pi(a)] \in \mathcal{L}^p(\mathcal{H})$ for all $a \in \mathcal{A}$, where  $\mathcal{L}^p
(\mathcal{H})$ is the \emph{Schatten ideal} of 
$p$-summable compact operators \cite{con2}. Fredholm modules should be thought of as 
representing
$K$-homology classes defined by abstract elliptic partial differential operators on the 
noncommutative space
$\mathcal{A}$.  

Spectral triples provide a refinement of Fredholm modules. Going from Fredholm modules 
to spectral triples is
akin to going from the \emph{conformal class} of a Riemannian metric to the
metric itself. Spectral triples simultaneously  provide a notion of
{\it Dirac operator} in noncommutative geometry, as well as a Riemannian
type {\it distance function} for noncommutative spaces as we shall explain next.

To motivate the definition of a spectral triple, we recall that the Dirac operator 
$D$ on a
compact Riemannian $\text{spin}^c$  manifold  acts  as an unbounded
selfadjoint operator on  the Hilbert space $L^2(M,S)$ of $L^2$-spinors on the manifold $M$. If we  let
$C^{\infty}(M)$ act  on $L^2(M,S)$ by multiplication operators, then  one can check that for any smooth
function $f$, the commutator $[D,f]=Df-fD$ 
extends to a bounded
operator on  $L^2(M,S)$. Now  the  \emph{geodesic distance} $d$ on $M$ can be recovered from the following 
beautiful {\it distance formula} of Connes \cite{acbook}:
\[d(p,q) = \textnormal{Sup} \{ |f(p)-f(q)|; \quad \parallel [D,f] \parallel \leq 1\}, \quad
 \,\,\, \forall p,q \in M.
\nonumber \]
The triple $(C^{\infty}(M), L^2(M,S), D)$ is a commutative example of
a spectral triple. Its general definition, in the odd case, is as
follows.
\begin{definition}
Let $\A$ be a unital algebra. An odd spectral triple on $\A$ is a triple $(\A,
\mathcal{H}, D)$ consisting  of a Hilbert space $\mathcal{H}$, a selfadjoint unbounded operator
$D: \text{Dom} (D) \subset \mathcal{H} \to \mathcal{H}$ with
compact resolvent, i.e.,
$(D+\lambda)^{-1} \in \mathcal{K}(\mathcal{H})  \, \text{for all} \,
\lambda \notin \mathbb{R},$
 and a unital 
representation  $\pi: \A \to \mathcal{L}(\mathcal{H})$  of $\A$ such that for
all $a\in \A$, the commutator $[D, \pi (a)]$ is defined on $\text{Dom} (D)$ and extends to  a bounded 
operator on $\mathcal{H}$.
\end{definition}

The finite summability assumption  for  Fredholm modules has a finer analogue
for spectral triples. For simplicity we shall assume  that $D$ is
invertible (in general, since $\text{Ker} (D)$ is finite dimensional, by
restricting to its orthogonal complement we can always reduce to this case).
  A  spectral
triple is called {\it finitely summable}  if for some $n\geq 1$,
\begin{equation}\label{stfs} |D|^{-n} \in \mathcal{L}^{1, \, \infty}(\mathcal{H}).
\end{equation}

A simple example of an odd spectral triple is
$(C^{\infty} (S^1), L^2 (S^1), D)$, where $S^1= \mathbb{R}/ 2 \pi \mathbb{Z}$ is the circle 
and $D$ is the unique selfadjoint extension of the  operator
$\frac{1}{i} \frac{d}{dx}$. The eigenvalues of $|D|$ are $|n|; n\in \mathbb{Z}$ which shows
that if we restrict $D$ to the orthogonal compliment of constant
functions then $|D|^{-1} \in \mathcal{L}^{1, \, \infty}(L^2(S^1)).$

Given a spectral triple $(\A, \mathcal{H}, D)$, one can obtain a
Fredholm module $(\A, \mathcal{H}, F)$ by choosing $F =\text{Sign}\,
(D)= D|D|^{-1}$. Connes' Hochschild character formula gives a local
expression for the Hochschild class of the Connes-Chern character of
$(\A, \mathcal{H}, F)$ in terms of $D$ itself. For this one has to assume that the spectral  triple 
$(\A, \mathcal{H}, D)$ is
 {\it regular}  in the sense that  for all $a\in \A$,
$$ a \, \, \text{and} \,\, [D, a] \, \in  \cap \,\text{Dom} (\delta^k)$$
where
 the derivation $\delta$ is  given by $\delta (x)= [|D|, x].$

 Now,  assuming \eqref{stfs} holds, Connes
defines an $(n+1)$-linear functional $\varphi$ on $\A$ by
$$\varphi (a^0, a^1, \dots, a^n)= \text{Tr}_{\omega} (a^0[D, a^1]\cdots
[D, a^n]|D|^{-n}). $$
It can be shown that $\varphi$ is a Hochschild $n$-cocycle on $\A$.
We recall that a Hochschild $n$-cycle $c \in Z_n (A, A)$ is an element
$c  =\sum a^0\otimes a^1 \otimes \dots \otimes a^n \in A^{\otimes (n+1)}$ such
that its Hochschild boundary $b (c)=0$, where $b$ is defined by \eqref{Hochschildboundary}.
The
 following result, known as {\it Connes' Hochschild character formula}, computes the Hochschild
 class of the Chern character by a local formula, i.e., in terms of
 $\varphi$:
\begin{theorem} 
\textnormal{\cite{acbook}}
Let  $(\A, \mathcal{H}, D)$ be a regular spectral
triple.
Let  $ F= \text{Sign}\, (D)$ denote the sign of $D$ and $\tau_n \in HC^n(\A)$  denote the  
Connes-Chern character of
$(\mathcal{H}, F)$. For every $n$-dimensional Hochschild cycle $c =\sum a^0\otimes a^1 
\otimes \dots \otimes a^n \in Z_n
(\A, \A)$, one has
$$\langle \tau_n, c\rangle =  \sum \varphi (a^0, a^1, \dots, a^n).$$
\end{theorem}

Identifying  the full cyclic cohomology class of the Connes-Chern character of 
$(\A, \mathcal{H}, D)$   by a
local formula is the content
of Connes-Moscovici's local index formula.
For this we have to assume the spectral triple satisfies  another
technical condition.
Let $\mathcal{B}$ denote the subalgebra of
$\mathcal{L}(\mathcal{H})$ generated by operators $\delta^k (a)$ and
$\delta^k ([D, a]), \, k\geq 1.$ A spectral triple is said to have a discrete 
{\it dimension spectrum} $\Sigma$
if $\Sigma \subset \mathbb{C}$ is discrete and for any $b\in
\mathcal{B}$ the function
$$ \zeta_b (z)= \text{Trace}(b |D|^{-z}), \, \,  \quad \text{Re}\, z > n,$$
extends to a holomorphic function on $\mathbb{C}\setminus  \Sigma$. It is further
assumed that $\Sigma$ is \emph{simple} in the sense that  $\zeta_b (z)$ has only
simple poles in $\Sigma$.

The local index formula of Connes and Moscovici \cite{conmos1} is given by the
following theorem (we have used the formulation in \cite{con4}):
\begin{theorem}
\textnormal{\cite{con4}}
1. The equality
$$ \int\!\!\!\!\!\!-P= \textnormal{Res}_{z=0} \, \textnormal{Trace} (P|D|^{-z})$$
defines a trace on the algebra generated by $\A$, $[D, \A]$, and $|D|^z,
\, z \in \mathbb{C}.$\\
2. There are  only a finite number of non-zero terms in the following
formula which defines the odd components $(\varphi_n)_{n=1, 3, \dots}$ of an odd 
cyclic cocycle in the $(b, B)$ bicomplex of
$\A$:
For each odd integer $n$ let
$$ \varphi_n (a^0, \cdots, a^n):= \sum_{k}c_{n, k} \int\!\!\!\!\!\!- a^0[D, a^1]^{(k_1)}\dots
[D, a^n]^{(k_n)}|D|^{-n-2|k|}$$
where $T^{(k)}:=  \nabla^k$ and $\nabla (T)= D^2T-TD^2$, $k$ is a
multi-index,
$|k|=k_1+\cdots + k_n$ and
$$ c_{n, k}:= (-1)^{|k|}\sqrt{2 i} (k_1!\dots k_n!)^{-1}((k_1+1)\cdots
(k_1+k_2+\cdots k_n))^{-1} \Gamma (|k|+\frac{n}{2}).$$
3. The pairing of the cyclic cohomology class  $(\varphi_n) \in HC^{*}(\A)$ with $K_1 (\A)$ 
gives the Fredholm index
of $D$ with coefficients in $K_1 (\A)$.
 \end{theorem}

Given an $n$-summable regular spectral triple $(\mathcal{A}, \mathcal{H}, D)$, the linear functional 
\[a \mapsto \text{Tr}_{\omega}(a |D|^{-n}) \nonumber \] 
defines a trace on the algebra $\mathcal{A}$ ($cf.$ Proposition \ref{twistedhypertrace}). 
Thus to deal with type III algebras
   which carry  no non-trivial traces,  the notion of spectral triple
   must be modified.  In \cite{conmos2} Connes and Moscovici define a
   notion of
{\it twisted spectral triple}, where the twist is afforded by an algebra automorphism 
(related to the modular automorphism group). More precisely,
one postulates that there exists an
automorphism $\sigma$ of $\mathcal{A}$ such that the twisted commutators
$Da-\sigma(a)D$
are bounded operators for all $a \in \mathcal{A}$. Here is the full definition:
\newtheorem{first}{Definition}[section]
\begin{first}
Let $\mathcal{A}$ be an algebra which is represented by bounded operators in a Hilbert space
$\mathcal{H}$, and $D$ be an unbounded selfadjoint operator in $\mathcal{H}$ with compact
resolvent. With $\sigma$ being an automorphism of $\mathcal{A}$, the triple $(\mathcal{A},
\mathcal{H}, D)$ is called a twisted spectral triple or a $\sigma$-spectral triple
if for any $a \in \mathcal{A}$, the twisted commutator
\[ [D,a]_{\sigma}:= Da-\sigma(a)D \nonumber\]
is defined on the domain of $D$, and extends to a bounded operator in $\mathcal{H}$.\\
A twisted spectral triple is said to be Lipschitz-regular if the twisted commutators
$|D|a-\sigma(a)|D|$ are bounded as well for all $a \in \mathcal{A}$.\\
A graded twisted spectral triple is one that is endowed with a grading operator
$\gamma=\gamma^* \in \mathcal{L}(\mathcal{H})$ such that $\gamma^2=id$, $\gamma$ commutes
with the action of $\mathcal{A}$, and anticommutes with $D$.
\end{first}
When the algebra $\mathcal{A}$ is involutive, the representation is assumed to be involutive
as well. In this case, it is natural to impose the following compatibility condition between
the automorphism and the involution:
\[\sigma(a)^*=\sigma^{-1}(a^*), \,\,\, \forall a \in \mathcal{A}.\nonumber\]

It is shown in \cite{conmos2} that in the twisted case, the Dixmier trace
induces a twisted trace on the algebra $\mathcal{A}$, but surprisingly, under some regularity conditions,
the Connes-Chern character of the phase space lands in ordinary cyclic cohomology. Thus its pairing
with ordinary $K$-theory makes sense, and it can be recovered as the index of Fredholm operators.
This suggests  the significance of developing
a local index formula for twisted spectral triples, $i.e.$ finding a formula for a cocycle, cohomologous to
the Connes-Chern character in the $(b,B)$-bicomplex, which is given in terms of twisted commutators and
residue functionals. We believe that this new theme of twisted spectral
triples, and type III noncommutative geometry in general,
shall dominate
the subject in near future.

For example, very recently a local index formula has been proved for a class of twisted spectral
 triples by Henri Moscovici \cite{mos1}.
We will discuss this result in detail in Section \ref{Local index formula and twisted spectral triples} 
of this paper.
This class is obtained by twisting an ordinary spectral triple $(\mathcal{A}, \mathcal{H}, D)$ by
a subgroup $G$ of \emph{scaling automorphisms} of the triple, $i.e.$ the set of all unitary operators 
$U \in \mathcal{U}
(\mathcal{H})$ such that $U \mathcal{A} U^*= \mathcal{A}$, and $UDU^*=\mu(U)D$, with $\mu(U) >0$. 
It is shown
that the crossed product algebra $\mathcal{A} \rtimes G$ admits an automorphism $\sigma$, given 
by the formula
$\sigma(aU)=\mu(U)^{-1}aU$, for all $a\in \mathcal{A}, U\in G$; and $(\mathcal{A} \rtimes G, \mathcal{H}, D)$
is a twisted spectral triple.

\section{Examples of twisted spectral triples} \label{Examples of twisted spectral triples}

In this section, we recall general methods to construct twisted spectral triples. 

\subsection{Perturbing spectral triples}

In \cite{conmos2}, it is shown that starting from an ordinary spectral triple
$(\mathcal{A}, \mathcal{H}, D)$ and a selfadjoint element $h=h^* \in \mathcal{A}$, the
perturbed triple
\[(\mathcal{A}, \mathcal{H}, D'), \,\,\, D'=e^hDe^h \nonumber \]
is a $\sigma$-spectral triple, where the automorphism $\sigma$ is given by
\[\sigma(a)=e^{2h}ae^{-2h}, \,\,\, a \in \mathcal{A}.\nonumber\]
In fact for any $a \in \mathcal{A}$, one has
\[D'a-\sigma(a)D'= e^hDe^ha-e^{2h}ae^{-2h}e^hDe^h=e^h[D, e^hae^{-h}]e^h.\nonumber\]
Therefore the twisted commutators $[D',a]_{\sigma}$ are bounded since the
commutators $[D, b]$ are bounded operators for all $b \in \mathcal{A}$. A concrete example of
this construction is obtained when one compares the Dirac operators of conformally equivalent
Riemannian metrics as follows:

\newtheorem{perturbex}[first]{Example} \label{perturbex}
\begin{perturbex} 
\rm{
Let $(M, g)$ be a compact Riemannian spin manifold and $D = 
D^g$ 
be the associated Dirac operator
acting on the Hilbert space of  $L^2$-spinors $\mathcal{H}= L^2(M, S^g)$. Let $g' = e^{-4h} g$ be a
conformally equivalent metric where $h \in C^{\infty}(M)$ is a selfadjoint element. It can be shown that
the gauge transform operator ${^g \hspace{-2pt}D^{g'}} = \beta_g^{g'} \circ D^{g'} \circ \beta^g_{g'}$ has 
the form 
\[  {^g \hspace{-2pt}D^{g'}} = e^{(n+1)h} \circ D^g \circ e^{-(n-1)h}. \nonumber \] 
After the canonical identification of the space of $g$-spinors with $g'$-spinors, it can be seen \cite{conmos2} 
that the gauge transform spectral triple is obtained from the original one by replacing $D$ with 
\[D' = e^hDe^h. \nonumber \]

}
\end{perturbex}

More generally, one can start from a $\sigma$-spectral triple
$(\mathcal{A}, \mathcal{H}, D)$ and a selfadjoint element $h=h^* \in \mathcal{A}$,
and investigate the properties of the perturbed triple
\[(\mathcal{A}, \mathcal{H}, D'), \,\,\, D'=e^hDe^h.\nonumber\]

\newtheorem{perturb}[first]{Lemma}
\begin{perturb}
Let $(\mathcal{A}, \mathcal{H}, D)$ be a $\sigma$-spectral triple, and $h=h^* \in \mathcal{A}$.
Then the perturbed triple
\[(\mathcal{A}, \mathcal{H}, D'), \,\,\, D'=e^hDe^h \nonumber\]
is a $\sigma'$-spectral triple where $\sigma' \in Aut(\mathcal{A})$ is given by
\[ \sigma'(a)=e^h\sigma(e^hae^{-h})e^{-h}, \,\,\, a \in \mathcal{A}. \nonumber\]
\begin{proof}
For any $a \in \mathcal{A}$, one has
\[D'a-\sigma'(a)D'= e^hDe^ha-e^h\sigma(e^hae^{-h})e^{-h}e^hDe^h=e^h[D, e^hae^{-h}]_{\sigma}e^h.
\nonumber\]
Therefore the twisted commutators $[D',a]_{\sigma'}$ are bounded since the twisted 
 commutators $[D, b]_{\sigma}$ are bounded operators for all $b \in \mathcal{A}$.
\end{proof}
\end{perturb}

\subsection{Twisted connections on crossed product algebras}

In \cite{fatkha1}, a method for constructing twisted spectral triples on crossed product algebras
is suggested. This method
uses 1-cocycles to construct automorphisms for  crossed product algebras, and explains how
one can obtain twisted traces, twisted derivatives, and  twisted connections on these
algebras. We shall explain this method and show that it can reconstruct an example of twisted 
spectral
triples, first given by Connes and Moscovici in \cite{conmos2}. First we recall some definitions:

\newtheorem{twisteddef}[first]{Definition}
\begin{twisteddef}\label{twisteddef}
Let $\sigma : \mathcal{A} \to \mathcal{A}$  be an automorphism of an algebra $\mathcal{A}$.
\begin{enumerate}
\item A $\sigma$-derivation on $\mathcal{A}$ is a linear map $\delta : \mathcal{A} \to \mathcal{A}$ 
such that
\begin{equation}
\delta (ab) = \delta (a) b + \sigma (a) \delta(b), \qquad \forall a,b \in \mathcal{A}. \nonumber
\end{equation}
\item A $\sigma$-trace on $\mathcal{A}$ is a linear map $\tau :\mathcal{A} \to \mathbb{C}$ such that
\begin{equation}
\tau (ab) = \tau (\sigma (b)a), \qquad \forall a,b \in \mathcal{A}. \nonumber
\end{equation}
\end{enumerate}
\end{twisteddef}

Now, let $A$ be an algebra with a right action of a group $\Gamma$ by
automorphisms:
\[ A \times \Gamma \to A,  \qquad
 (a,\gamma) \mapsto a \cdot \gamma . \nonumber\]
We consider the algebraic crossed product $A \rtimes \Gamma$ with the standard multiplication:
\[ (a \otimes \gamma)(b \otimes \mu) = (a \cdot \mu) b \otimes \gamma \mu, \quad \,\,\, a,b \in A, 
\,\,\, \gamma, \mu \in \Gamma.\nonumber \]
Let $Z(A)$ denote the center of the algebra $A$, and $A^*$ its group of invertible elements.

\newtheorem{chain}[first]{Definition}
\begin{chain}\label{chain}
\begin{enumerate}
\item A map $j : \Gamma \to Z(A) \cap A^*$ is a 1-cocycle if
\begin{equation} \label{eq:jmap}
j(\gamma \mu) = \big ( (j\gamma) \cdot \mu \big ) j \mu  , \quad \,\,\, \forall \gamma, \mu \in \Gamma .
\nonumber
\end{equation}
\item Given a map $j : \Gamma \to A$, a linear functional $\tau :A \to \mathbb{C}$ is said to have the 
change of variable property with respect to $j$, if
\[ \tau((a \cdot \gamma)j \gamma)= \tau (a),  \quad \,\,\, \forall a \in A, \,\,\, \gamma \in \Gamma.\nonumber\]
\end{enumerate}
\end{chain}
Notice that \eqref{eq:jmap} amounts to saying that $j$ is a
(multiplicative) group $1-$cocycle for $H^1(\Gamma , Z(A) \cap A^*).$

\newtheorem{automor}[first]{Proposition}
\begin{automor}\label{twistedderivation}
\textnormal{\cite{fatkha1}}
Let $A$ be an algebra with a right action of a group $\Gamma$ by automorphisms, and $j :
 \Gamma \to Z(A) \cap A^*$ be a 1-cocycle.
\begin{enumerate}
\item The map $\sigma : A \rtimes \Gamma \to A \rtimes \Gamma$  given by
\[ \sigma (a \otimes \gamma) = \big ( (j \gamma^{-1})\cdot \gamma \big ) a \otimes \gamma \nonumber\]
is an automorphism.
\item Let $\delta : A \to A$ be a derivation such that $\delta (a \cdot \gamma) = (\delta (a) \cdot \gamma) \,
 j \gamma$ for all $a \in A, \gamma \in \Gamma$. Then for any $s=1,2, \dots$, the map $\delta'_s :
  A \rtimes \Gamma \to A \rtimes \Gamma$ defined by
\[ \delta'_s (a \otimes \gamma) = \Big ( \delta \Big ( (a \cdot \gamma^{-1}) (j \gamma^{-1})^s \Big ) 
(j \gamma^{-1})^{-s} \otimes 1 \Big ) \Big ( 1 \otimes \gamma \Big ) \nonumber\]
is a $\sigma$-derivation on $A \rtimes \Gamma$. Also $\delta'_s \circ \sigma = \sigma \circ \delta'_s$ 
if $\delta \circ j = 0$.
\item If $\tau : A \to \mathbb{C}$ is a trace such that $\tau \circ \delta = 0$, then $\tau' \circ \delta'_s =0$ 
where the linear functional $\tau' : A \rtimes \Gamma \to \mathbb{C}$ is defined by
\[
\tau' (a \otimes \gamma) = 0 \,\, if \,\,  \gamma \neq 1, \,\,\, and \,\,\, \tau' (a \otimes 1) = \tau(a).\nonumber\]
Also, if $\tau$ has the change of variable property with respect to $j$, then $\tau'$ is a $\sigma$-trace 
on $A \rtimes \Gamma.$

\end{enumerate}

\end{automor}

\newtheorem{jmapex}[first]{Example}
\begin{jmapex}
\begin{enumerate}
\item Let $M$ be a smooth oriented manifold and $\omega$ a volume form on $M$. Let $\Gamma = 
Diff(M)$ 
be the group of diffeomorphisms of $M$. The map $j: \Gamma \to C^{\infty} (M)$ defined by
\[ \gamma ^* (\omega) = j(\gamma) \omega \nonumber\]
is easily seen to be a 1-cocycle.
\item Let $\chi : \Gamma \to \mathbb{C}^*$ be a 1-dimensional character of a group $\Gamma$ which 
acts on an algebra $A$ by automorphisms, and let $j(\gamma) = \chi (\gamma) 1_{A}$. Then $j$ is a 
1-cocycle,
and a derivation $\delta : A \to A$ is compatible with $j$ if and only if
$\delta (a \cdot \gamma) = \chi (\gamma) \,\, \delta(a) \cdot \gamma$ for
any $a \in A$, and $\gamma \in \Gamma.$
\end{enumerate}
\end{jmapex}

\newtheorem{connectiondef}[first]{Definition}
\begin{connectiondef}\label{connectiondef}
Let $\mathcal{A}$ be an algebra, $\delta : \mathcal{A} \to \mathcal{A}$ a derivation, and $E$ a left 
$\mathcal{A}$-module. A linear map $\nabla : E \to E$ is said to be a connection if it satisfies the 
Leibniz rule, i.e.

\[ \nabla (a  \xi) = \delta (a)  \xi + a  \nabla (\xi), \,\,\, \forall a \in \mathcal{A}, \,\,\, \xi \in E.\nonumber \]

If $\sigma : \mathcal{A} \to \mathcal{A}$ is an automorphism and $\delta : \mathcal{A} \to \mathcal{A}$  a 
$\sigma$-derivation, then a linear map $\nabla : E \to E$ is said to be a twisted connection if it satisfies the
 twisted Leibniz rule:
\[ \nabla (a  \xi) = \delta (a) \xi + \sigma (a)  \nabla (\xi), \,\,\, \forall a \in \mathcal{A}, \,\,\, \xi \in E.\nonumber \]

\end{connectiondef}
This notion of twisted connection comes from \cite{pol}. \\

Given a twisted connection $\nabla$, one can try to define a twisted  spectral triple by letting $D = \nabla$. 
Then $\nabla a - \sigma (a) \nabla = \delta (a)$ shows the boundedness condition
is satisfied provided $\delta (a)$ acts by a bounded operator. \\

In the following proposition, $A$ is an algebra endowed with a right action of a
group $\Gamma$ by automorphisms, with representations $\pi : A \to \text{End} \,(E)$, and
$\rho: \Gamma \to \text{GL} (E)$,  defining  a covariant system, i.e.
\[\pi (a \cdot \gamma ) = \rho (\gamma^{-1})  \pi(a)  \rho(\gamma), \,\,\, \forall a \in A, \,\,\, \gamma \in 
\Gamma. \nonumber \]
Then  we obtain  a representation $\pi' : A \rtimes \Gamma \to \text{End} (E)$ given by
\[\pi' (a \otimes \gamma) = \rho(\gamma) \pi (a), \,\,\,  \forall a \in A, \,\,\, \gamma \in \Gamma.\nonumber\]
Also let  $j : \Gamma \to Z(A) \cap A^*$ be a 1-cocycle and
$\delta : A \to A$  a derivation  as in Proposition \ref{twistedderivation}.
Therefore we have an automorphism $\sigma : A \rtimes \Gamma \to A \rtimes
\Gamma $ and we fix a $\sigma$-derivation $\delta'_s : A \rtimes \Gamma \to A \rtimes \Gamma$ for 
some $s \in \mathbb{N}$.

\newtheorem{twistedconnection}[first]{Proposition}
\begin{twistedconnection}\label{twistedconnection}
\textnormal{\cite{fatkha1}}
A connection  $\nabla: E \to E $ for $A$ is a twisted connection for $A \rtimes \Gamma$ if and only if
\begin{equation} \label{twistedconncond}
\nabla  \rho (\gamma) = \rho (\gamma) \Big ( (s-1) \, \pi \big (\delta (j \gamma^{-1}\cdot \gamma) \big ) + 
\nabla \pi (j \gamma^{-1}\cdot \gamma)  \Big ), \,\,\, \forall \gamma \in \Gamma.
\end{equation}
\end{twistedconnection}

\newtheorem{reals}[first]{Remark}
\begin{reals}
\rm{
Let $s > 0$ be a real number and assume $\big ( j (\gamma) \big )^s \in A$ is defined for all 
$\gamma \in \Gamma$. Propositions \ref{twistedderivation} and \ref{twistedconnection} 
continue to hold for these values of $s$ as well. For this we need the extra condition 
$\delta (x^s) = s x^{s-1} \delta (x)$ to hold for all $x = j(\gamma) $, $\gamma \in \Gamma$.}
\end{reals}

\newtheorem{foliationex}[first]{Example}
\begin{foliationex}
\rm{ Let $ C^{\infty} (S^1)$ be the algebra of smooth functions on the circle 
$S^1 = \mathbb{R} / \mathbb{Z}$, and $\Gamma \subset Diff^+ (S^1)$ a group 
of orientation preserving diffeomorphisms of the circle as in \cite{conmos2}.  
We represent the algebra $C^{\infty} (S^1)$ by bounded operators in the 
Hilbert space $L^2(S^1)$ by
\[ (\pi (g) \, \xi) (x) = g(x) \, \xi (x), \,\,\, \forall g \in C^{\infty} (S^1), \,\,\,  \xi \in 
L^2(S^1), \,\,\, x \in \mathbb{R} / \mathbb{Z}. \nonumber\]
Define a representation of $\Gamma$ by bounded operators in  $L^2(S^1)$ by
\[ (\rho(\phi^{-1}) \, \xi )(x) = \phi'(x)^{\frac{1}{2}} \, \xi(\phi(x)), \,\,\, \forall \phi \in 
\Gamma, \,\,\,  \xi \in L^2(S^1), \,\,\, x \in \mathbb{R} / \mathbb{Z}. \nonumber\]
The group $\Gamma$ acts on $ C^{\infty} (S^1)$ from right by composition and one can
easily check  that the above representations give a covariant system which yields the
representation of $ C^{\infty} (S^1) \rtimes \Gamma$ as  in  \cite{conmos2}.
The map $j : \Gamma \to C^{\infty} (S^1)$ defined by $j(\phi) = \phi'$ is a 1-cocycle 
and the derivation $\delta: C^{\infty} (S^1)
\to C^{\infty} (S^1), \delta (f) = \frac{1}{i} f'$ is compatible with $j$. Now by using
Proposition \ref{twistedderivation}, we obtain an automorphism $\sigma$
of $ C^{\infty} (S^1) \rtimes \Gamma$ which agrees with the automorphism
in \cite{conmos2}, and a twisted derivation $\delta'_{\frac{1}{2}}$. Note that since it 
is possible to take the square root of the elements in the image of $j$ in this example, 
we can let $s=\frac{1}{2}$ to obtain a twisted derivation. Now if we let $\nabla = 
\frac{1}{i} \frac{d}{dx}$, one can see that the equality
\eqref{twistedconncond} holds, therefore $\nabla$ is a twisted connection for 
$C^{\infty}(S^1) \rtimes \Gamma$  by Proposition \ref{twistedconnection}.

Also the linear map $\tau : C^{\infty}(S^1) \to \mathbb{C}$ defined by
\[ \tau (g) = \int_{\mathbb{R} / \mathbb{Z}} g(x) \, dx, \,\,\, \forall g \in C^{\infty}(S^1), \nonumber\]
is a trace which has the change of variable property with respect to $j$, and $\tau \circ \delta = 0$.
 Therefore by Proposition \ref{twistedderivation}, one obtains a twisted trace $\tau' : C^{\infty}(S^1) 
 \rtimes \Gamma \to \mathbb{C} $ such that $\tau' \circ \delta'_{\frac{1}{2}} = 0.$

}
\end{foliationex}

\subsection{Twisting spectral triples by scaling automorphisms} 
\label{Twisting spectral triples by scaling automorphisms}

The last example of the preceding subsection gives a special case of a class of twisted spectral
triples that arise naturally in conformal geometry \cite{mos1}. Let $(M, g)$ be a connected compact 
Riemannian spin manifold of dimension $n$ and $D = D_g$ 
be the associated Dirac operator acting on the Hilbert space of  $L^2$-spinors $\mathcal{H}= L^2(M, S^g)$. 
Let $SCO(M, [g])$ denote the Lie group of diffeomorphisms of $M$ that preserve the conformal structure $[g]$
consisting of all Riemannian metrics that are conformally equivalent to $g$, the orientation, and the spin
structure, and let $G = SCO(M, [g])_0$ denote the connected component of the identity. In \cite{mos1}, using
a suitable automorphism of the crossed product algebra $C^{\infty}(M) \rtimes G$, a twisted spectral triple of
the form $(C^{\infty}(M) \rtimes G,   L^2(M, S^g), D)$ is constructed. Similarly, by endowing
$\mathbb{R}^n$ with the Euclidean metric, and considering the group $G$ of conformal transformations of 
$\mathbb{R}^n$, a twisted
spectral triple  is constructed over the crossed product algebra $C_c^{\infty}(\mathbb{R}^n) \rtimes G$. An abstract
formulation of this class of twisted spectral triples leads in \cite{mos1} to the idea of twisting an ordinary spectral 
triple by its \emph{scaling automorphisms} or  \emph{conformal similarities} which we explain in this subsection.  

Using scaling automorphisms of a spectral triple $(\mathcal{A}, \mathcal{H}, D)$, 
one can
construct a twisted spectral triple \cite{mos1}. The set of scaling automorphisms, also called 
conformal similarities, 
of a spectral triple $(\mathcal{A}, \mathcal{H}, D)$, denoted by $\textnormal{Sim}(\mathcal{A}, 
\mathcal{H}, D)$,
consists of all unitary operators $U$ on $\mathcal{H}$ such that
\[U \mathcal{A} U^*= \mathcal{A}, \,\,\, \textnormal{and} \,\,\, UDU^*=\mu(U)D, \,\,\,
\textnormal{for some} \,\,\,  \mu(U) >0. \nonumber \]
It is easy to see that $\textnormal{Sim}(\mathcal{A}, \mathcal{H}, D)$ is a group and the map $\mu :
\textnormal{Sim}(\mathcal{A}, \mathcal{H}, D) \to (0, \infty)$ is a character. We fix a subgroup
$G \subset \textnormal{Sim}(\mathcal{A}, \mathcal{H}, D)$ and let $\mathcal{A}_G$ be the crossed product
algebra $\mathcal{A} \rtimes G$. It is shown in \cite{mos1} that the formula
\[\sigma(aU)= \mu(U)^{-1}aU, \,\,\, \forall a\in \mathcal{A}, U \in G, \nonumber \]
defines an automorphism of $\mathcal{A}_G$, and $(\mathcal{A}_G, \mathcal{H}, D)$ is a twisted
spectral triple. In fact the twisted commutators
\[ [D,aU]_{\sigma}:=DaU - \sigma(aU) D =[D, a ]U \nonumber \]
are bounded operators for all $a \in \mathcal{A}, U \in G$.

For this class of twisted spectral triples,
one can form the crossed product algebra  $\Psi(\mathcal{A}\rtimes G, \mathcal{H}, D) := \Psi(\mathcal{A},
\mathcal{H}, D) \rtimes G$, where $\Psi(\mathcal{A}, \mathcal{H}, D)$ is the algebra of pseudodifferential
operators associated to the base spectral triple \cite{conmos1, hig1, hig2}, and under the \emph{extended 
simple dimension
spectrum hypothesis}, the residue functional ${\int_{D} \hspace{-17pt} -} \,:\Psi(\mathcal{A}_G,
\mathcal{ H}, D) \to \mathbb{C}$ given by
\[{\int_{D} \hspace{-15pt} -}P:= \textnormal{Res}_{z=0} \textnormal{Trace} (P |D|^{-2z}) \nonumber \]
is a trace \cite{mos1}. This seems to be in agreement with the twisted analogue of the Adler-Manin trace
\cite{adl, man, fatkha1}: For a triple $(A, \sigma, \delta)$ consisting of an algebra $A$, an
algebra automorphism $\sigma:A \to A$, and a $\sigma$-derivation $\delta :A\to A$, we define
an algebra of formal twisted pseudodifferential symbols $\Psi(A, \sigma, \delta)$ whose elements
are formal series of the form
\[ \sum_{n=-\infty}^N a_n \xi^n, \,\,\, \,\,\, N \in \mathbb{Z}, \, a_n \in A, \, \forall n \leq N. \nonumber  \]
The multiplication in this algebra is essentially derived from these relations
\[ \xi a = \sigma(a)\xi + \delta(a), \forall a \in A, \,\,\, \xi \xi^{-1}= \xi^{-1} \xi =1. \nonumber  \]
In \cite{fatkha1}, we prove that starting from a $\delta$-invariant twisted trace on $A$, the
noncommutative residue is a trace on $\Psi(A, \sigma, \delta)$:

\newtheorem{nctrace}[first]{Theorem}
\begin{nctrace} \label{nctrace}
\textnormal{\cite{fatkha1}}
Let $\sigma : A \to A$ be an automorphism of an algebra $A$. If $\tau: A \to \mathbb{C}$ 
is a $\sigma$-trace  and $\delta: A \to A$ is a $\sigma$-derivation such that
$\tau \circ \delta = 0$, then the linear functional $\textnormal{Res}:
\Psi (A, \sigma , \delta)  \to \mathbb{C}$ defined by
\[ \textnormal{Res} \, \big ( \sum_{i=- \infty}^{n} a_i \xi^i \big ) =
\tau (a_{-1}) \nonumber \]
is a trace on $\Psi (A, \sigma , \delta)$.
\end{nctrace}

\section{Properties of twisted spectral triples} \label{Properties of twisted spectral triples}

In this section, we recall the basic properties of twisted spectral triples \cite{conmos2}.

\subsection{Twisted trace}

Given a  $\sigma$-spectral triple $(\mathcal{A}, \mathcal{H}, D)$ with
$D^{-1}\in \mathcal{L}^{n, \infty}$, it is observed in \cite{conmos2}
that
\[ D^{-n} a - \sigma^{-n}(a) D^{-n} \in \mathcal{L}_0^{1, \infty}(\mathcal{H}), \,\,\, \forall a \in
\mathcal{A}. \nonumber \]
Here $\mathcal{L}_0^{1, \infty}(\mathcal{H})$ is the ideal
\[\{  T \in \mathcal{K}(\mathcal{H}) ; \quad \sum_{i=0}^{N} \mu_i(T) = o(\log N) \}, \nonumber \]
on which the Dixmier trace $\textnormal{Tr}_{\omega}$ vanishes. Then it follows that the Dixmier trace
induces a twisted hypertrace on $\mathcal{A}$:

\newtheorem{twistedhypertrace}[first]{Proposition}
\begin{twistedhypertrace} \label{twistedhypertrace} \textnormal{\cite{conmos2}}
Let $(\mathcal{A}, \mathcal{H}, D)$ be a $\sigma$-spectral triple  with
$D^{-1}\in \mathcal{L}^{n, \infty}$. Then the linear functional $\varphi: \mathcal{A} \to
\mathbb{C}$ defined by
\[\varphi(a) = \textnormal{Tr}_{\omega}(aD^{-n}), \,\,\, a \in \mathcal{A}  \nonumber \]
is a $\sigma^{n}$-hypertrace, $i.e.$
\[ \textnormal{Tr}_{\omega} (T a D^{-n}) = \textnormal{Tr}_{\omega} (\sigma^{n}(a) T D^{-n}) \nonumber \]
for any $a \in \mathcal{A}, T\in \mathcal{L}(\mathcal{H}).$ In particular $\varphi$ is a twisted trace on
$\mathcal{A}.$ Note that for Lipschitz-regular twisted spectral triples, the same holds when $D^{-n}$
is replaced by $|D|^{-n}$.
\end{twistedhypertrace}

\subsection{Connes-Chern character}

Let $(\mathcal{A}, \mathcal{H}, D)$ be a Lipschitz-regular $\sigma$-spectral triple such that
$D^{-1} \in \mathcal{L}^{n, \infty}(\mathcal{H})$ for some $n \in \mathbb{N}$. Here
\[ \mathcal{L}^{n, \infty}(\mathcal{H}) = \{  T \in \mathcal{K}(\mathcal{H}) ; \quad
\sum_{i=0}^{N} \mu_i(T) = O(N^{1-\frac{1}{n}}) \}, \,\,\, \text{if} \,\,\, n>1, \nonumber \]
\[ \mathcal{L}^{1, \infty}(\mathcal{H}) = \{  T \in \mathcal{K}(\mathcal{H}) ; \quad
\sum_{i=0}^{N} \mu_i(T) = O(\log N) \},\nonumber \]
where for any compact operator $T \in \mathcal{K}(\mathcal{H})$, its singular values are
written in decreasing order: $\mu_1(T) \geq \mu_2(T) \geq \cdots \geq 0.$

We recall from \cite{conmos2} that passage to the phase $F= D|D|^{-1}$ of such a twisted spectral triple
gives a finitely summable Fredholm module which has a well-defined Connes-Chern character
in cyclic cohomology given by
\begin{equation} \label{Connes-Chern}
\Phi_F(a_0, a_1, \dots, a_n) = \textnormal{Trace}(\gamma F [F,a_0][F,a_1] \cdots [F,a_n]), \,\,\, 
a_0, a_1, \dots , a_n \in \mathcal{A}, \nonumber
\end{equation}
where $\gamma=id$ in the ungraded case.
Since for any $a \in \mathcal{A}$
\[ [F, a]=|D|^{-1}([D,a]_{\sigma} - [|D|, a]_{\sigma} F ), \nonumber \]
the commutators $[F, a_i]$ are differentials of the same order as $D^{-1}$ and the operator
\[\gamma F [F,a_0][F,a_1] \cdots [F,a_n] \nonumber \] 
is a trace class operator. One can see
that $\Phi_F$ is a cyclic cocycle, $i.e.$
\[ \Phi_F (a_n, a_0, \dots, a_{n-1}) = (-1)^n \Phi_F (a_0, a_1, \dots, a_n), \,\,\, \forall
a_0, a_1, \dots , a_n \in \mathcal{A},  \nonumber \]
\[b \Phi_F = 0,  \nonumber \]
where $b$ is the Hochschild coboundary operator:
\begin{eqnarray} \label{Hochschildcoboundary}
b \Phi_F(a_0, a_1, \dots, a_{n+1}) &=& \sum_{i=0}^{n} (-1)^i \Phi_F(a_0, \dots, a_i
a_{i+1}, \dots, a_{n+1})
+ \nonumber \\
&&(-1)^{n+1}\Phi_F(a_{n+1}a_0, a_1, \dots, a_n), \nonumber 
\end{eqnarray}
for all $a_0, a_1, \dots , a_{n+1} \in \mathcal{A}$.

Moreover, if $D^{-1} \in \mathcal{L}^{n, \infty}(\mathcal{H})$
for an even $n\in \mathbb{N}$, then the Connes-Chern
character can be defined without the Lipschitz-regularity assumption:

\newtheorem{Connes-ChernwnLR}[first]{Proposition}
\begin{Connes-ChernwnLR} \label{Connes-ChernwnLR}
\textnormal{\cite{conmos2}}
Let $(\mathcal{A}, \mathcal{H}, D)$ be a graded $\sigma$-spectral triple with
$D^{-1} \in \mathcal{L}^{n, \infty}(\mathcal{H})$ for some even $n \in \mathbb{N}$. Then the
multilinear functional $\Phi_{D, \sigma}$ defined by
\begin{equation} \label{ConnesChernenLRF}
\Phi_{D, \sigma}(a_0, a_1, \dots, a_n)=\textnormal{Trace}(\gamma D^{-1} [D, a_0]_{\sigma} 
D^{-1} [D, a_1]_{\sigma}
\cdots D^{-1} [D, a_n]_{\sigma}) 
\end{equation}
for $a_0, \dots, a_n\in \mathcal{A}$, is a cyclic cocycle.
\end{Connes-ChernwnLR}

\subsection{Index pairing and $K$-theory}
In \cite{conmos2}, it is shown that in the twisted case, the index pairing with ordinary
$K$-theory makes sense and it is given by the pairing of the Connes-Chern character with $K$-theory.

Given a  graded $\sigma$-spectral triple $(\mathcal{A}, \mathcal{H}, D)$  with
$D^{-1} \in \mathcal{L}^{n, \infty}(\mathcal{H})$ for some even $n \in \mathbb{N}$, one can write an
orthogonal decomposition for the Hilbert space $\mathcal{H}$ using the grading $\gamma$:
\[ \mathcal{H}=\mathcal{H}_+ \oplus \mathcal{H}_- , \nonumber \]
where
\[ \mathcal{H}_+ = \{ x \in \mathcal{H}; \quad \gamma(x)=x\},  \,\,\, \mathcal{H}_- =
\{ x \in \mathcal{H}; \quad \gamma(x)=-x\}. \nonumber \]
Now a close look at the Connes-Chern character $\Phi_{D, \sigma}$ defined by 
\eqref{ConnesChernenLRF}
shows the existence of a pair of Fredholm modules over $\mathcal{A}$ and two cocycles as follows.
With respect to the decomposition $\mathcal{H}=\mathcal{H}_+ \oplus \mathcal{H}_-$, we can write
\begin{equation}
D =
\begin{bmatrix}
0 & D_- \\
D_+ & 0 
\end{bmatrix},
\qquad
a =
\begin{bmatrix}
a_+ & 0 \\
0 & a_-
\end{bmatrix},
\,\,\, \forall a \in \mathcal{A}, \nonumber
\end{equation}
since we have $D\gamma = - \gamma D$, and $a \gamma = \gamma a$, for all $a \in \mathcal{A}$.
Then define two Hilbert spaces
\[ \widetilde{\mathcal{H}}^+ = \mathcal{H}_+ \oplus \mathcal{H}_+, \,\,\,
\widetilde{\mathcal{H}}^- = \mathcal{H}_- \oplus \mathcal{H}_-.  \nonumber \]
 There are representations $\pi^+$ and $\pi^-$ of the algebra $\mathcal{A}$ on
 $\mathcal{L}({\widetilde{\mathcal{H}}}^+)$ and $\mathcal{L}({\widetilde{\mathcal{H}}}^-)$ respectively given by
\begin{equation}
\pi^{\pm}(a)=
\begin{bmatrix}
a_{\pm} & 0 \\
0 & D_{\pm}^{-1} \sigma(a)_{\mp} D_{\pm}
\end{bmatrix}. \nonumber
\end{equation}
In \cite{conmos2}, letting
\[ F = \begin{bmatrix}
0 & I_{\pm} \\
I_{\pm} & 0
\end{bmatrix}, \nonumber \]
where $I_+$ and $I_-$ are identity operators on $\mathcal{H}_+$ and $\mathcal{H}_-$ 
respectively, it is shown
that the commutators $[F^{\pm}, \pi^{\pm}(a)]$ are compact operators, hence they obtain a pair of 
Fredholm modules
over the algebra $\mathcal{A}$. It is also shown that for any idempotent $e \in \mathcal{A}$, the 
bounded closure of
operators $D_{\pm}^{-1}\sigma(e)_{\mp}D_{\pm}$ denoted by $f_{\pm}$ are idempotents, and
$f_{\pm} e_{\pm}: e_{\pm} \mathcal{H}_{\pm} \to f_{\pm} \mathcal{H}_{\pm}$  are Fredholm operators. 
Since the index
depends only on the $K$-theory class of the idempotent, they define a pair of index
maps by
\begin{equation}
\textnormal{Index}^{\pm}[e] = \textnormal{Index} (f_{\pm}e_{\pm}), \nonumber
\end{equation}
for all idempotents $e \in M_N(\mathcal{A})$.

On the other hand, the cyclic cocycle $\Phi_{D, \sigma}$ defined by \eqref{ConnesChernenLRF} 
consists of a
pair of cocycles $\Phi^{\pm}_{D, \sigma}$ given by
\begin{eqnarray}
&& \Phi^{\pm}_{D, \sigma}(a^0, \dots, a^n) \nonumber \\ 
&=& \textnormal{Trace}( D_{\pm}^{-1} (D_{\pm}a^0_{\pm} - \sigma(a^0)_{\mp} D_{\pm})
\cdots D_{\pm}^{-1} (D_{\pm}a^n_{\pm} - \sigma(a^n)_{\mp} D_{\pm})), \nonumber
\end{eqnarray}
for $a^0, \dots, a^n \in \mathcal{A}.$

The following proposition states that the index pairing can be expressed as the pairing of these 
cocycles with
$K$-theory:

\newtheorem{PairingwithK-theroy}[first]{Proposition}
\begin{PairingwithK-theroy} \textnormal{\cite{conmos2}}
Given a graded $\sigma$-spectral triple $(\mathcal{A}, \mathcal{H}, D)$  with
$D^{-1} \in \mathcal{L}^{n, \infty}(\mathcal{H})$ for some even $n \in \mathbb{N}$, and an
idempotent $e \in M_N(\mathcal{A})$, one has
\[ \textnormal{Index}^{\pm}[e] = \Phi^{\pm}_{D, \sigma}(e, \dots , e). \nonumber \]
If $e^*= \sigma(e)$, then
\[ \textnormal{Index}^+[e] = - \textnormal{Index}^-[e]. \nonumber \]
\end{PairingwithK-theroy}

If there exists
a strongly continuous 1-parameter group of isometric automorphisms $\{ \sigma_t \}_{t \in \mathbb{R}} $ 
with an analytic
extension coinciding with $\sigma$ at $t=-i$, then the above index maps coincide. This assumption is 
denoted by (1PG).
Such an analytic extension defined on a dense subalgebra $\mathcal{O}$ of $\mathcal{A}$ is assured to 
exist by a theorem of
Bost in \cite{bos}.

\newtheorem{indexmaps}[first]{Proposition}
\begin{indexmaps} \textnormal{\cite{conmos2}}
If $\mathcal{A}$ satisfies (1PG), then
\[ \textnormal{Index}^+[e] = - \textnormal{Index}^-[e], \,\,\, \forall e \in M_N(\mathcal{A}). \nonumber \]
\end{indexmaps}

Accordingly, in \cite{conmos2}, the relation between the cyclic cocycles $\Phi^{\pm}_{D, \sigma}$ has 
been studied
under the (1PG) assumption.

\newtheorem{cocycles}[first]{Theorem}
\begin{cocycles} \textnormal{\cite{conmos2}}
Let $(\mathcal{A}, \mathcal{H}, D)$ be a graded $\sigma$-spectral triple as above and assume that 
$\mathcal{A}$ satisfies
(1PG). Then
\[ [\Phi^{-}_{D, \sigma}] = - [(\Phi^{+}_{D, \sigma})^*] \in HP^{ev}(\mathcal{O}), \nonumber \]
where
\[(\Phi^{+}_{D, \sigma})^*(a_0, \dots, a_n) := \overline{ \Phi^{+}_{D, \sigma}(a^*_n, \dots, a^*_0)}, 
\,\,\, \forall a_0,
\dots, a_n \in \mathcal{O}. \nonumber \]
\end{cocycles}

We note that in the proof of the latter \cite{conmos2}, the homotopy invariance of the Connes-Chern character of
a finitely summable Fredholm module, established in Lemma 1, in Part I, section 5 of  \cite{con2} 
plays a crucial role.

\subsection{Local Hochschild cocycle} \label{Local Hochschild cocycle}
In \cite{conmos2}, as a first step to extend the local index formula \cite{conmos1, hig1, hig2}
to twisted spectral triples, using the Dixmier trace, a
local Hochschild cocycle is constructed for twisted spectral triples:

\newtheorem{Hochschildcocycle}[first]{Proposition}
\begin{Hochschildcocycle} \label{Hochschildcocycle}
\textnormal{\cite{conmos2}}
Let $(\mathcal{A}, \mathcal{H}, D)$ be a graded $\sigma$-spectral triple such that
$D^{-1} \in \mathcal{L}^{n, \infty}(\mathcal{H})$ for some even $n \in \mathbb{N}$.
Then the $n+1$-linear form on $\mathcal{A}$ defined by
\begin{eqnarray} \label{Hochschild}
&& \Psi_{D,\sigma}(a_0, a_1, \dots ,a_n) \nonumber \\
&=&    \textnormal{Tr}_{\omega} ( \gamma a_0 [D, \sigma^{-1}(a_1)]_{\sigma} \cdots
[D, \sigma^{-n}(a_n)]_{\sigma} |D|^{-n}) \nonumber
\end{eqnarray}
for $a_0, \dots, a_n \in \mathcal{A}$, is a Hochschild cocycle.

In the ungraded case, for a Lipschitz-regular $\sigma$-spectral triple of odd summability
degree, the Hochschild cocycle is given by
\begin{eqnarray}
&& \Psi_{D,\sigma}(a_0, a_1, \dots ,a_n) \nonumber \\
&=&    \textnormal{Tr}_{\omega} (  a_0 [D, \sigma^{-1}(a_1)]_{\sigma} \cdots
[D, \sigma^{-n}(a_n)]_{\sigma} |D|^{-n}), \nonumber
\end{eqnarray}
for any $a_0, \dots, a_n \in \mathcal{A}$.
\end{Hochschildcocycle}

The above cocycle identities are proved in \cite{conmos2} using Proposition \ref{twistedhypertrace}
and the fact that for any $a,b \in \mathcal{A}$:
\[ [D, ab]_{\sigma}=[D, a]_{\sigma}b+ \sigma(a)[D, b]_{\sigma}. \nonumber \]

In \cite{conmos2, mos1}, the above local Hochschild cocycle is obtained in a heuristic manner as follows. 
For an
ordinary spectral triple $(\mathcal{A}, \mathcal{H}, D)$, consider the local Hochschild cocycle
\[ \Psi_D(a_0, \dots, a_n) = \textnormal{Tr}_{\omega}(\gamma a_0 [D, a_1] \cdots [D, a_n] D^{-n}), 
\nonumber \]
for any $a_0, \dots, a_n \in \mathcal{A}$.
One can move $D^{-n}$ to the left and distribute it among the factors to write this cocycle in the form:
\[ \Psi_D(a_0, \dots, a_n) = \textnormal{Tr}_{\omega}(\gamma a_0 (Da_1D^{-1}-a_1) \cdots 
(D^n a_nD^{-n} - D^{n-1}a_nD^{-n+1}) D^{-n}). \nonumber\]
In order to construct a local Hochschild cocycle for twisted spectral triples, they replace each $D^kaD^{-k}$ 
in the
latter by $D^k\sigma^{-k}(a)D^{-k}$ and reverse the process of moving the $D^{-n}$ to the left which leads to 
the above local
Hochschild cocycle for twisted spectral triples. We note that the Connes-Chern character introduced in 
Proposition 
\ref{Connes-ChernwnLR} can be obtained in a similar manner.

\section{Hochschild class of the Connes-Chern character} \label{Hochschild class of the Connes-Chern character}

For any algebra $\mathcal{A}$, there is an obvious pairing between the space of
Hochschild $n$-cochains $C^n(\mathcal{A}, \mathcal{A}^*)$, $i.e.$ the space of
$(n+1)$-linear functionals on $\mathcal{A}$, and the space of Hochschild $n$-chains
$\mathcal{A}^{\otimes (n+1)}$, given by 
\[ \langle \varphi, a_0 \otimes a_1 \otimes \cdots \otimes a_n \rangle
:= \varphi (a_0, \dots , a_n). \nonumber \]
This pairing satisfies 
\[ \langle b \varphi, c \rangle = \langle  \varphi, b (c) \rangle \nonumber \]
where the Hochschild operators for the cochains and chains are given by:
\begin{eqnarray} \label{Hochschildcoboundary}
b \varphi(a_0, a_1, \dots, a_{n+1}) &=& \sum_{i=0}^{n} (-1)^i \varphi(a_0, \dots, a_i
a_{i+1}, \dots, a_{n+1})
+ \nonumber \\
&&(-1)^{n+1}\varphi(a_{n+1}a_0, a_1, \dots, a_n), \nonumber
\end{eqnarray}
\begin{eqnarray} \label{Hochschildboundary}
b (a_0 \otimes a_1 \otimes \dots \otimes a_{n+1}) &=& \sum_{i=0}^{n} 
(-1)^i a_0 \otimes \dots \otimes a_i
a_{i+1} \otimes \dots \otimes a_{n+1}
+ \nonumber \\
&&(-1)^{n+1} a_{n+1}a_0 \otimes a_1 \otimes \dots \otimes a_n,
\end{eqnarray}
for $a_0, \dots, a_{n+1} \in \mathcal{A}$.
It follows that, if two cocycles $\varphi_1$ and $\varphi_2$ are cohomologous with
$\varphi_1 - \varphi_2 = b \psi$, then they yield the same value on any Hochschild
cycle since we have
\[ \langle \varphi_1, c \rangle  - \langle \varphi_2, c \rangle
= \langle b\psi, c \rangle = \langle \psi, b(c) \rangle = 0 \,\,\, \text{if} \,\,\, b(c)=0. \nonumber \]

As we saw in Section \ref{Properties of twisted spectral triples}, 
to any twisted spectral 
triple, in particular to an ordinary spectral triple, with certain conditions 
one can associate the Connes-Chern character and a 
local Hochschild cocycle. Connes' character formula (or Connes' Hochschild character theorem) states 
that in the case of ordinary spectral triples, these two cocycles yield the same value
on any Hochschild cycle \cite{acbook, gravarfig, hig2}.     

In this section the analogue of Connes' character formula for 
the class of twisted spectral triples introduced in Subsection \ref{Twisting spectral triples by scaling 
automorphisms} is investigated which we explain next following \cite{fat1}. 

\subsection{Connes' character formula and twisting by scaling automorphisms}

In Subsection \ref{Twisting spectral triples by scaling automorphisms}, we saw that using conformal 
similarities of a spectral triple
$(\mathcal{A}, \mathcal{H}, D)$, one can
construct a twisted spectral triple ($cf.$ \cite{mos1}). We recall that the set of conformal
similarities of a spectral
triple $(\mathcal{A}, \mathcal{H}, D)$, denoted by $\textnormal{Sim}(\mathcal{A}, \mathcal{H}, D)$,
consists of all unitary operators $U$ on $\mathcal{H}$, such that
\[U \mathcal{A} U^*= \mathcal{A}, \,\,\, \textnormal{and} \,\,\, UDU^*=\mu(U)D, \,\,\,
\textnormal{for some} \,\,\,  \mu(U) >0. \nonumber \]
Recall that $\textnormal{Sim}(\mathcal{A}, \mathcal{H}, D)$ is a group and the map $\mu :
\textnormal{Sim}(\mathcal{A}, \mathcal{H}, D) \to (0, \infty)$ is a character. We fix a subgroup
$G \subset \textnormal{Sim}(\mathcal{A}, \mathcal{H}, D)$ and let $\mathcal{A}_G$ be the crossed product
algebra $\mathcal{A} \rtimes G$. The formula
\[\sigma(aU)= \mu(U)^{-1}aU, \,\,\, \forall a\in \mathcal{A}, U \in G, \nonumber\]
defines an automorphism of $\mathcal{A}_G$, and $(\mathcal{A}_G, \mathcal{H}, D)$ is a
$\sigma$-spectral triple \cite{mos1}.
We will assume that the base spectral triple $(\mathcal{A}, \mathcal{H}, D)$ is regular, $i.e.$
the operators $\mathcal{A}, [D, \mathcal{A}]$ are in the domain of all powers of the derivation
$\delta(\cdot)=[|D|, \cdot]$. We also assume that $D^{-1}\in \mathcal{L}^{n,\infty}(\mathcal{H})$
for some fixed even number $n \in 2 \mathbb{N}$, and that the
twisted spectral triple  $(\mathcal{A}_G, \mathcal{H}, D)$ is graded.

From the regularity of the base spectral triple, it easily follows that the twisted spectral
triple  $(\mathcal{A}_G, \mathcal{H}, D)$ is Lipschitz-regular. Therefore, by passage to the phase, 
one can associate the
Connes-Chern character to the Fredholm module $(\mathcal{A}_G, \mathcal{H}, F=D|D|^{-1})$. 
To recall, this is
a cyclic $n$-cocycle given by
\[ \Phi_F(a_0U_0, a_1U_1, \dots, a_nU_n) = \textnormal{Trace}(\gamma F [F,a_0U_0][F,a_1U_1] 
\cdots [F,a_nU_n]) \nonumber \]
for all $a_i \in \mathcal{A}, U_i \in G, i = 0, \dots , n $. Here $\gamma$ denotes the grading.
Also there is a Hochschild $n$-cocycle  given by
\begin{eqnarray} \label{Hoch}
&& \Psi_{D,\sigma}(a_0 U_0, a_1 U_1, \dots ,a_n U_n) \nonumber \\
&=&    \textnormal{Tr}_{\omega} ( \gamma a_0 U_0 [D, \sigma^{-1}(a_1 U_1)]_{\sigma} \cdots
[D, \sigma^{-n}(a_n U_n)]_{\sigma} |D|^{-n}) \nonumber
\end{eqnarray}
for all $a_i \in \mathcal{A}, U_i \in G, i = 0, \dots ,n $.

We recall from \cite{gravarfig} how
\[ \Phi_F(a_0U_0, a_1U_1, \dots, a_nU_n) = \textnormal{Trace}(\gamma F [F,a_0U_0][F,a_1U_1]
\cdots[F,a_nU_n]) \nonumber \]
can be approximated by the trace of finite rank operators using a cutoff. Let $g:\mathbb{R}
\to \mathbb{R}$ be a smooth function which takes the value $1$ on $[0,\frac{1}{2}]$, decreases to
$0$ on $[\frac{1}{2},1]$, is $0$ for $t>1$, $g(-t) = g(t)$ for $t<0$, and $\int_{0}^{\infty}g'(u) \text{d}u=-1$.
Define the operators $A_t = g(t|D|) $ for $t >0$. The operators $A_t$, $t >0$, are positive, finite rank, and
satisfy $P_{1/2t} \leq A_t \leq P_{1/t}$, where $P_N$ is the spectral
projector of $|D|$ on the interval $[0,N]$.  Given a Hochschild $n$-cycle
\[ c= \sum_{j=1}^{k} a_{0j}U_{0j}\otimes a_{1j}U_{1j} \otimes \cdots \otimes a_{nj}U_{nj},\nonumber \]
it is proved in \cite{gravarfig} that
\begin{equation} \label{estimate}
 \Phi_F(c) = 2 \lim_{t \downarrow 0} \Psi_t (c), 
\end{equation}
where
\begin{eqnarray}
&&\Psi_t (a_0U_0, a_1U_1, \dots, a_nU_n) \nonumber \\
&:=& - \textnormal{Trace} (\gamma a_0U_0 [F, a_1U_1] \cdots [F, a_{n-1}U_{n-1}]
F[A_t,a_nU_n]). \nonumber
\end{eqnarray}
The proof of this is essentially based on two facts ($cf.$ \cite{gravarfig}). The first is that the
operator trace is normal, and the second is that, since $b(c)=0$, one has
\begin{eqnarray}
&& \sum_{j=1}^{k} a_{0j}U_{0j}[F, a_{1j}U_{1j}] \cdots [F, a_{n-1,j}U_{n-1,j}] a_{nj}U_{nj}  \nonumber  \\
&=&  \sum_{j=1}^{k}a_{nj}U_{nj} a_{0j}U_{0j}[F, a_{1j}U_{1j}] \cdots [F, a_{n-1,j}U_{n-1,j}]. \nonumber
\end{eqnarray}

For convenience,  we drop the index $j$ and the summation in the formula for the Hochschild cycle $c$,
and simply write
\[c = a_{0}U_{0}\otimes a_{1}U_{1} \otimes \cdots \otimes a_{n}U_{n}. \nonumber \]

\newtheorem{first1}[first]{Lemma}
\begin{first1}\label{approximate1}
\textnormal{\cite{fat1}}
The operator $a_0U_0[F, a_1U_1] \cdots [F, a_{n-1}U_{n-1}]F|D|^{n-1}$ is bounded.
\begin{proof}
Using the identity $U[F,a] = [F,UaU^*]U$, for all $a\in A, U \in G$, one can see that
\begin{eqnarray}
&&a_0U_0[F, a_1U_1] \cdots [F, a_{n-1}U_{n-1}]F|D|^{n-1} \nonumber \\
&=&   \mu(U_0U_1 \cdots U_{n-1})^{n-1}a_0 [F, U_0a_1U^*_0] \cdots \nonumber \\
&&  [F,U_0U_1 \cdots U_{n-2} a_{n-1}U^*_{n-2} \cdots U^*_1U^*_0 ] F |D|^{n-1}U_0U_1
\cdots U_{n-1}.  \nonumber
\end{eqnarray}
The boundedness of the latter follows from the fact that the operator
\[a'_0[F, a'_1] \cdots [F, a'_{n-1}]F|D|^{n-1} \nonumber \]
is bounded for any $a'_0, \dots , a'_{n-1} \in A$, which is proved in \cite{gravarfig}.
\end{proof}
\end{first1}

For convenience let
\[R= -\gamma a_0U_0[F, a_1U_1] \cdots [F, a_{n-1}U_{n-1}]F|D|^{n-1} \in \mathcal{L}(H). \nonumber \]
Then we have:
\[ \Psi_t (a_0U_0, a_1U_1, \dots, a_nU_n) =  \textnormal{Trace} (R|D|^{-(n-1)}[A_t,a_nU_n]). \nonumber  \]

In the sequel we will impose the following extra condition on the Hochschild cycle $c  = 
a_{0}U_{0}\otimes a_{1}U_{1} \otimes \cdots \otimes a_{n}U_{n}$ (note that we have dropped the summation):

\newtheorem{limit}[first]{Condition}
\begin{limit}\label{condition}
We shall assume that
\[\lim_{t \downarrow 0} \textnormal{Trace} (R|D|^{-(n-1)}[A_t,a_nU_n]) = \lim_{t \downarrow 0} 
\textnormal{Trace}(R|D|^{-(n-1)}
[A_t,a_n]U_n). \nonumber \]
For example, if $\mu(U_n)=1$ then this condition is satisfied.
\end{limit}

The function $\Psi'_t (a_0U_0, a_1U_1, \dots, a_nU_n) := \textnormal{Trace} (R|D|^{-(n-1)}[A_t,a_n]U_n)$ is
continuous on the interval $\frac{1}{t} \geq e$, and has a limit as $\frac{1}{t} \to \infty$,
therefore it is bounded. Hence the evaluation of the state $\omega$ on the corresponding
element in the $C^*$-algebra $C_b([e, \infty))/C_0([e, \infty))$ will yield this limit which will
be denoted by 
\[\lim_{t^{-1} \to \omega} \Psi'_t(a_0U_0, a_1U_1, \dots, a_nU_n). \nonumber \] 
To compute this limit, one can use the following proposition of Connes. For a detailed discussion,
we refer the reader to \cite{gravarfig}.

\newtheorem{computelimit}[first]{Proposition}
\begin{computelimit}\label{computelimit}
Let $f:[0,\infty) \to [0,\infty)$ be a continuous function, $p>1$, and $\sum_k m_k(f)e^{pk} < \infty$,
where for each $k$
\[m_k(f)= \sup \{f(u); \quad k \leq \log u \leq k+1\}. \nonumber \]
Then $M_p = \int_{0}^{\infty} f(u) u^{p-1} \textnormal{d}u$ is finite and
\[ \lim_{t^{-1} \to \omega} t^p \textnormal{Trace} (f(t|D|)S) = M_p 
\textnormal{Tr}_{\omega} (S |D|^{-p}), \nonumber \]
provided $S \in \mathcal{L}(\mathcal{H}) $, and $D^{-1} \in \mathcal{L}^{p, \infty}(\mathcal{H}).$
\end{computelimit}

Also we use Lemma 10.29 of \cite{gravarfig}. 
\newtheorem{cutofflemma}[first]{Lemma}
\begin{cutofflemma}
If $g(t)= h(t)^2$ where $h \in \mathcal{D}(\mathbb{R})$ is also a cutoff, and if
$a \in \mathcal{A}$, then
\[ \parallel  [g(t|D|), a] - \frac{1}{2} \{g'(t|D|), t \delta{a} \} \parallel_{n-} = o(t)
\,\,\, \textnormal{as} \,\,\, t \downarrow 0.\nonumber \]
\end{cutofflemma}
\noindent We note that $\parallel \cdot \parallel_{n-}$ is defined in \cite{gravarfig}.

The argument following the above lemma in \cite{gravarfig} applies verbatim to our case by
simply replacing $R$ by $U_nR$, and it follows that
\begin{eqnarray} \label{relating}
 && \lim_{t^{-1} \to \omega} \Psi'_t(a_0U_0, a_1U_1, \dots, a_nU_n) \nonumber \\
 &=& n \textnormal{Tr}_{\omega}( \gamma a_0U_0 [F,a_1U_1] \cdots [F,a_{n-1}U_{n-1}] D^{-1} 
 \delta(a_n) U_n ) \nonumber \\
&=&  n \textnormal{Tr}_{\omega} (\gamma a_0U_0 [F,a_1U_1] \cdots [F,a_{n-1}U_{n-1}] 
 \delta_{\sigma}\sigma^{-1}
(a_nU_n) D^{-1}), 
\end{eqnarray}
where $\delta(x)=[|D|, x]$ and $\delta_{\sigma}(x) = [|D|, x]_{\sigma}$.

The last step is to prove that the local Hochschild cocycle $\Psi_{D,\sigma}$ is cohomologous to
$\zeta^{\sigma}_n$ defined by
\begin{eqnarray}
&& \zeta^{\sigma}_n(a_0U_0, a_1U_1, \dots, a_nU_n) \nonumber \\
&=&  n \textnormal{Tr}_{\omega} (\gamma a_0U_0 [F,a_1U_1] \cdots [F,a_{n-1}U_{n-1}]  
\delta_{\sigma}\sigma^{-1}
(a_nU_n) D^{-1}) \nonumber
\end{eqnarray}
for any $a_i \in \mathcal{A},$ and $U_i \in G$, $i = 0, \dots, n.$ This will be achieved by defining the 
cochains $\zeta^{\sigma}_1, \dots ,\zeta^{\sigma}_{n-1}$ by defining
$\zeta^{\sigma}_k(a_0U_0, a_1U_1, \dots, a_nU_n)$ as
\[n \textnormal{Tr}_{\omega} (\gamma a_0U_0 [F,a_1U_1] \cdots [F,a_{k-1}U_{k-1}] D^{-1} 
\delta_{\sigma}
(a_kU_k)[F,a_{k+1}U_{k+1}] \cdots [F,a_{n}U_{n}]) \nonumber \]
for any $a_i \in \mathcal{A}, U_i \in G$, $i=0, \dots, n$. Then using the following two lemmas, 
one can see that $\Psi_{D,\sigma}$ is cohomologous to $\zeta^{\sigma}_n$, hence they yield
the same value on any Hochschild cycle.

\newtheorem{ncocycles0}[first]{Lemma}
\begin{ncocycles0} \label{ncocycles0}
\textnormal{\cite{fat1}}
The cochains $\zeta^{\sigma}_{1}, \dots ,\zeta^{\sigma}_{n}$ are mutually cohomologous Hochschild
cocycles.
\begin{proof}
First we show that $b\zeta^{\sigma}_n=0.$ For any $a_i \in \mathcal{A}, U_i \in G,$ we have
\begin{eqnarray}
&& b\zeta^{\sigma}_n (a_0 U_0, \dots ,a_{n+1} U_{n+1}) \nonumber \\
&=&   \textnormal{Tr}_{\omega}(\gamma a_0U_0 a_1U_1 [F,a_2U_2]\cdots [F,a_nU_n][|D|,
\sigma^{-1}(a_{n+1}U_{n+1})]_{\sigma}D^{-1}) -\nonumber \\
&&  \textnormal{Tr}_{\omega}(\gamma a_0U_0 [F,a_1U_1a_2U_2]\cdots [F,a_nU_n][|D|,
\sigma^{-1}(a_{n+1}U_{n+1})]_{\sigma}D^{-1}) -\nonumber \\
&&   \cdots + \nonumber \\
&&  \textnormal{Tr}_{\omega}(\gamma a_0U_0 [F,a_1U_1]\cdots [F,a_{n-1}U_{n-1}][|D|,
\sigma^{-1}(a_nU_na_{n+1}U_{n+1})]_{\sigma}D^{-1}) \nonumber \\
&&   -\textnormal{Tr}_{\omega}(\gamma a_{n+1}U_{n+1}a_0U_0 [F,a_1U_1]\cdots [F,a_{n-1}U_{n-1}][|D|,
\sigma^{-1}(a_nU_n)]_{\sigma}D^{-1}) \nonumber 
\end{eqnarray}
\begin{eqnarray}
&=& \textnormal{Tr}_{\omega}(\gamma a_0U_0 a_1U_1 [F,a_2U_2]\cdots [F,a_nU_n][|D|,
\sigma^{-1}(a_{n+1}U_{n+1})]_{\sigma}D^{-1}) \nonumber \\
&&  - \textnormal{Tr}_{\omega}(\gamma a_0U_0 a_1U_1 [F,a_2U_2]\cdots [F,a_nU_n][|D|,
\sigma^{-1}(a_{n+1}U_{n+1})]_{\sigma}D^{-1}) \nonumber \\
&&  - \textnormal{Tr}_{\omega}(\gamma a_0U_0 [F,a_1U_1]a_2U_2 \cdots [F,a_nU_n][|D|,
\sigma^{-1}(a_{n+1}U_{n+1})]_{\sigma}D^{-1}) \nonumber \\
&&  \cdots \nonumber \\
&&  +\textnormal{Tr}_{\omega}(\gamma a_0U_0 [F,a_1U_1]\cdots [F,a_{n-1}U_{n-1}]a_nU_n[|D|,
\sigma^{-1}(a_{n+1}U_{n+1})]_{\sigma}D^{-1})  \nonumber \\
&& + \textnormal{Tr}_{\omega}(\gamma a_0U_0 [F,a_1U_1]\cdots [F,a_{n-1}U_{n-1}][|D|,
\sigma^{-1}(a_nU_n)]_{\sigma} \nonumber \\ 
&& \qquad \qquad \sigma^{-1}(a_{n+1}U_{n+1})D^{-1}) \nonumber \\
&&  -\textnormal{Tr}_{\omega}(\gamma a_{n+1}U_{n+1}a_0U_0 [F,a_1U_1]\cdots [F,a_{n-1}U_{n-1}][|D|,
\sigma^{-1}(a_nU_n)]_{\sigma}D^{-1}). \nonumber
\end{eqnarray}
The latter vanishes because of successive cancellations, where the last two terms cancel each
other since using the identity
\[ D^{-1}a_{n+1}U_{n+1} - \sigma^{-1}(a_{n+1}U_{n+1}) D^{-1} = -D^{-1}[D,
\sigma^{-1}(a_{n+1}U_{n+1})]_{\sigma}D^{-1}, \nonumber \]
one obtains a trace class operator which vanishes under the Dixmier trace.

Now we introduce cochains $\eta^{\sigma}_1, \dots, \eta^{\sigma}_{n-1}$ such that
$\zeta^{\sigma}_{k} - \zeta^{\sigma}_{k+1} = b \eta^{\sigma}_k$ for all $k=1, \dots , n-1$,
and this will finish the proof.

Using the identities
\[  D^{-1} \delta_{\sigma}(a_kU_k) -  \delta_{\sigma}\sigma^{-1}(a_kU_k) D^{-1} = - D^{-1}
\delta([D,a_k])D^{-1}U_k, \nonumber \]
\[ |D|^{-1}[F,a_jU_j]-[F,\sigma^{-1}(a_jU_j)]|D|^{-1} = -|D|^{-1}[F, \delta_{\sigma}
\sigma^{-1}(a_jU_j)]|D|^{-1},  \nonumber \]
\[F[F,a_jU_j]=-[F,a_jU_j]F, \nonumber \]
we can move $D^{-1}$ in the formula for $\zeta^{\sigma}_k$ to the right under the Dixmier trace
and obtain the following expression. \begin{eqnarray}
&&\zeta^{\sigma}_k(a_0U_0, a_1U_1, \dots, a_nU_n)  \nonumber \\
 &=&  n \textnormal{Tr}_{\omega} (\gamma a_0U_0 [F,a_1U_1] \cdots [F,a_{k-1}U_{k-1}] \nonumber \\
 && D^{-1} \delta_{\sigma}(a_kU_k)[F,a_{k+1}U_{k+1}] \cdots [F,a_{n}U_{n}])
\nonumber \\
&=&  (-1)^{n-k}n \textnormal{Tr}_{\omega}(\gamma a_0U_0 [F,a_1U_1] \cdots [F,a_{k-1}U_{k-1}] \nonumber \\
 &&   \delta_{\sigma}(\sigma^{-1}(a_kU_k))[F,\sigma^{-1}(a_{k+1}U_{k+1})]
\cdots [F,\sigma^{-1}(a_{n}U_{n})]D^{-1}  ) . \nonumber
\end{eqnarray}
Therefore $(\zeta^{\sigma}_k - \zeta^{\sigma}_{k+1})(a_0U_0, a_1U_1, \dots, a_nU_n) $ is equal to
\begin{eqnarray}
(-1)^{n-k}n \textnormal{Tr}_{\omega} (\gamma a_0U_0 [F, a_1U_1] \cdots [F, a_{k-1}U_{k-1}] R^{\sigma}_k
[F,\sigma^{-1}(a_{k+2}U_{k+2})]\nonumber \\
 \cdots  [F, \sigma^{-1}(a_nU_n)]D^{-1}), \nonumber
\end{eqnarray}
where 
\[R^{\sigma}_k= \delta_{\sigma}\sigma^{-1}(a_kU_k)[F, \sigma^{-1}(a_{k+1}U_{k+1})]+
[F,a_kU_k]\delta_{\sigma}\sigma^{-1}(a_{k+1}U_{k+1}). \nonumber \]
Now we define
\begin{eqnarray}
&& \eta^{\sigma}_k(a_0U_0, \dots ,a_{n-1}U_{n-1}) \nonumber \\
&=&  (-1)^{k} n \textnormal{Tr}_{\omega} (\gamma a_0U_0 [F, a_1U_1] \cdots [F, a_{k-1}U_{k-1}] \nonumber \\
&&   [F, \delta_{\sigma}\sigma^{-1}(a_kU_k)][F, \sigma^{-1}(a_{k+1}U_{k+1})]
\cdots [F, \sigma^{-1}(a_{n-1}U_{n-1})]D^{-1} ). \nonumber
\end{eqnarray}
Finally, using the identity
\begin{eqnarray}
&& [F, \delta_{\sigma}\sigma^{-1}(a_kU_ka_{k+1}U_{k+1})]  \nonumber \\
&=&  R^{\sigma}_k + [F, \delta_{\sigma}\sigma^{-1}(a_kU_k)] \sigma^{-1}(a_{k+1}U_{k+1}) +
a_k U_k [F, \delta_{\sigma}\sigma^{-1}(a_{k+1}U_{k+1})], \nonumber
\end{eqnarray}
one can see that $\zeta^{\sigma}_{k} - \zeta^{\sigma}_{k+1} = b \eta^{\sigma}_k.$
\end{proof}
\end{ncocycles0}

\newtheorem{ncocycles}[first]{Lemma}
\begin{ncocycles} \label{ncocycles}
\textnormal{\cite{fat1}}
The cochain $\Psi_{D,\sigma} - \frac{1}{n} (\zeta^{\sigma}_1 + \cdots + \zeta^{\sigma}_n)$ is a
Hochschild coboundary.
\begin{proof}
Let $a_0 U_0, a_1 U_1, \dots ,a_n U_n \in \mathcal{A}_G$.
Since
\begin{eqnarray}
&&|D|^{-1} [D, a_jU_j]_{\sigma} - [D, \sigma^{-1}(a_jU_j)]_{\sigma} |D|^{-1} \nonumber \\
&=& |D|^{-1} [D, a_j]U_j - \mu(U_j)[D, a_j] U_j |D|^{-1}  \nonumber \\
&=&   |D|^{-1} [D, a_j]U_j - [D, a_j]  |D|^{-1} U_j  \nonumber \\
&=&  - |D|^{-1} \delta ([D,a]) |D|^{-1} U, \nonumber
\end{eqnarray}
in the expression for $\Psi_{D,\sigma}$ we can replace each $[D, \sigma^{-1}(a_jU_j)]_{\sigma}
|D|^{-1}$ by \\
$|D|^{-1} [D, a_jU_j]_{\sigma}$. Therefore
\begin{eqnarray} \label{Hoch}
&& \Psi_{D,\sigma}(a_0 U_0, a_1 U_1, \dots ,a_n U_n) \nonumber \\
&=&    \textnormal{Tr}_{\omega} ( \gamma a_0 U_0 [D, \sigma^{-1}(a_1 U_1)]_{\sigma} \cdots
[D, \sigma^{-n}(a_n U_n)]_{\sigma} |D|^{-n}) \nonumber \\
&=&   \textnormal{Tr}_{\omega} ( \gamma a_0 U_0 [D, \sigma^{-1}(a_1 U_1)]_{\sigma} |D|^{-1} \cdots
[D, \sigma^{-1}(a_n U_n)]_{\sigma} |D|^{-1}).
\end{eqnarray}
Also we have
\[ [D, \sigma^{-1}(a_jU_j)]_\sigma |D|^{-1} = [F, a_jU_j ] + \delta_{\sigma} \sigma^{-1}(a_jU_j)
D^{-1} + [F, \delta_{\sigma} \sigma^{-1}(a_jU_j)]  |D|^{-1}. \nonumber  \]
Since
\begin{eqnarray}
[F, \delta_{\sigma} \sigma^{-1}(a_jU_j)] &=& \mu(U_j) [F, \delta(a_j) U_j] \nonumber \\
&=& \mu(U_j) [F, \delta(a_j)] U_j \nonumber \\
&=& \mu(U_j) |D|^{-1} (\delta([D,a_j]) - \delta^2(a_j)F)U_j, \nonumber
\end{eqnarray}
the terms containing $[F, \delta_{\sigma} \sigma^{-1}(a_jU_j)]  |D|^{-1}$ yield trace class
operators which vanish under the
Dixmier trace. So we can replace each 
\[ [D, \sigma^{-1}(a_jU_j)]_\sigma |D|^{-1} \nonumber  \] 
by
\[[F, a_jU_j ] + \delta_{\sigma} \sigma^{-1}(a_jU_j) D^{-1}. \nonumber \] 
Therefore \eqref{Hoch} is the sum of
$2^n$ terms. The term 
\[\textnormal{Tr}_{\omega} (\gamma a_0 U_0 [F,a_1U_1] \cdots [F,a_nU_n]) \nonumber \] 
is zero
because one can write $a_0U_0 = F[F,a_0U_0]+Fa_0U_0F$, and 
\[\gamma F[F,a_0U_0][F,a_1U_1] \cdots [F,a_nU_n] \nonumber \] 
is trace class. Therefore the term is equal to
\begin{eqnarray}
&& \textnormal{Tr}_{\omega} (\gamma Fa_0U_0F[F,a_1U_1] \cdots [F,a_nU_n]) \nonumber \\
&=&  (-1)^n \textnormal{Tr}_{\omega} (\gamma Fa_0U_0[F,a_1U_1] \cdots [F,a_nU_n]F) \nonumber \\
&=&   -\textnormal{Tr}_{\omega} (F \gamma a_0U_0[F,a_1U_1] \cdots [F,a_nU_n]F) \nonumber \\
&=&   -\textnormal{Tr}_{\omega} (\gamma a_0U_0[F,a_1U_1] \cdots [F,a_nU_n]). \nonumber
\end{eqnarray}
Hence it has to be zero.

The terms having exactly one factor of the form $\delta_{\sigma} \sigma^{-1} (a_jU_j) D^{-1}$ add
up to $\frac{1}{n} (\zeta^{\sigma}_1 + \cdots + \zeta^{\sigma}_n)$, and to finish the proof we show
that the terms with more than one factors of the form $\delta_{\sigma} \sigma^{-1} (a_jU_j) D^{-1}$
are Hochschild coboundaries. For example let us consider the case when two consecutive factors of the
above form yield the term
\begin{eqnarray} \label{ncocylcesid}
&& \textnormal{Tr}_{\omega}( \gamma a_0U_0 [F,a_1U_1] \cdots [F,a_{j-1}U_{j-1}] \delta_{\sigma} 
\sigma^{-1}(a_jU_j) \nonumber \\
&&  D^{-1} \delta_{\sigma} \sigma^{-1}(a_{j+1}U_{j+1})D^{-1}
[F,a_{j+2}U_{j+2}]\cdots [F,a_nU_n] ) \nonumber \\
&=&  \textnormal{Tr}_{\omega}(\gamma a_0U_0 [F,a_1U_1] \cdots [F,a_{j-1}U_{j-1}] \delta_{\sigma}
\sigma^{-1}(a_jU_j) \nonumber \\
&&  \delta_{\sigma} \sigma^{-2}(a_{j+1}U_{j+1}) D^{-2} [F,a_{j+2}U_{j+2}]
\cdots [F,a_nU_n] ) \nonumber \\
&=&  \textnormal{Tr}_{\omega}(\gamma a_0U_0 [F,a_1U_1] \cdots [F,a_{j-1}U_{j-1}] \delta_{\sigma}
\sigma^{-1}(a_jU_j) \nonumber \\
&&  \delta_{\sigma} \sigma^{-2}(a_{j+1}U_{j+1}) [F,\sigma^{-2}(a_{j+2}U_{j+2})]\cdots
[F,\sigma^{-2}(a_nU_n)] D^{-2}).
\end{eqnarray}
Now using the identity
\[\delta^2_{\sigma} \sigma^{-2}(a_j a_{j+1}) = \delta^2_{\sigma} \sigma^{-2}(a_j) \sigma^{-2}
(a_{j+1}) + 2 \delta_{\sigma} \sigma^{-1}(a_j) \delta_\sigma  \sigma^{-2} (a_{j+1}) + a_j
\delta^2_{\sigma} \sigma^{-2} (a_{j+1}), \nonumber \]
one can see that \eqref{ncocylcesid} is equal to $b\varphi^{\sigma}_j(a_0U_0, a_1U_1, \dots, a_nU_n)$
where
\begin{eqnarray}
&& \varphi^{\sigma}_j(a_0U_0, a_1U_1, \dots, a_{n-1}U_{n-1}) \nonumber \\
&=& \frac{(-1)^j}{2} \textnormal{Tr}_{\omega} (\gamma a_0U_0[F,a_1U_1] \cdots [F,a_{j-1}U_{j-1}] \nonumber \\
&& \delta^2_{\sigma}\sigma^{-2}(a_jU_j)[F,\sigma^{-2}(a_{j+1}U_{j+1})] \cdots
[F,\sigma^{-2}(a_{n-1}U_{n-1})]D^{-2}). \nonumber
\end{eqnarray}
\end{proof}

\end{ncocycles}

Hence, considering \eqref{estimate}, Condition \ref{condition}, and \eqref{relating} we have proved 
the following:
\newtheorem{mainthm}[first]{Theorem}
\begin{mainthm} \label{mainthm}
\textnormal{\cite{fat1}}
The cyclic cocycle $\Phi_F$ and the local Hochschild cocycle $2 \Psi_{D,\sigma}$ yield the same value
on any Hochschild $n$-cycle
\[c= \sum_{j=1}^{k} a_{0j}U_{0j}\otimes a_{1j}U_{1j} \otimes \cdots \otimes a_{nj}U_{nj} \nonumber \]
that satisfies Condition \ref{condition}.
\end{mainthm}

\section{A local index formula for twisted spectral triples} 
\label{Local index formula and twisted spectral triples}
In \cite{mos1}, an Ansatz for a local index formula for twisted spectral triples is given and
its validity for twisted spectral triples obtained from scaling automorphisms of spectral triple, 
has been verified. This formula is given in terms of residue functionals and twisted commutators. 
In this section we sketch very briefly some of the ideas in \cite{mos1}. 

\subsection{Moscovici's Ansatz}
In Subsection \ref{Local Hochschild cocycle}, we explained the heuristic manner for obtaining
the local Hochschild cocycle for twisted spectral triples from the one for ordinary spectral
triples. In a rather similar manner one can obtain Moscovici's Ansatz for a local index formula for
twisted spectral triples from the original local index formula of Connes and Moscovici discussed in 
Section
\ref{From spectral triples to twisted spectral triples}. The Ansatz suggests that for a $\sigma$-spectral
triple $(\mathcal{A}, \mathcal{H}, D)$, the twisted version of
the local character should be given by
\begin{eqnarray}
&& \varphi^{\sigma}_n (a_0, \dots, a_n) := \nonumber \\
&& \sum_{k}c_{n, k} \int\!\!\!\!\!\!- a_0[D, \sigma^{-2k_1-1}(a_1)]_{\sigma}^{(k_1)}\cdots
[D, \sigma^{-2(k_1+ \cdots + k_n) - n}(a_n)]_{\sigma}^{(k_n)}|D|^{-n-2|k|}, \nonumber
\end{eqnarray}
where the \emph{iterated twisted commutators} $[D, a]_{\sigma}^{(k)}$ are defined in \cite{mos1}.
Here again
\[\int\!\!\!\!\!\!- P = \textnormal{Res}_{z=0} \textnormal{Trace} (P |D|^{-2z}). \nonumber \]
Therefore an analogue of the simple dimension spectrum hypothesis is assumed. Namely, 
one needs to assume that
there exists a discrete subset of the complex plane such that for any operator $P$ in an algebra 
of twisted pseudodifferential 
operators associated to the twisted spectral triple $(\mathcal{A}, \mathcal{H}, D, \sigma)$, $P |D|^{-2z}$ 
is a trace class operator 
provided $\text{Re}(z)$ is large enough, and $ \zeta_P(z) := \textnormal{Trace}(P|D|^{-2z})$ has a meromorphic 
extension to the plane with 
at most simple poles in this discrete set. We note that, as it is emphasized in \cite{mos1}, there is \emph{no} 
canonical way of constructing an algebra of  twisted pseudodifferential 
operators for a general twisted spectral triple. Also, 
in the twisted case,  for formulating analogue of the regularity condition, one needs to postulate that the automorphism 
$\sigma$ has an extension to a larger algebra which contains  the \emph{twisted differential forms}.    
Another necessary condition for
the validity of the Ansatz is a $\sigma$-invariance of the residue functionals which is related to the Selberg
principle for reductive Lie groups \cite{mos1}. 

\subsection{A local index formula for spectral triples twisted by scaling automorphisms}
\label{}
Let $(\mathcal{A}_G, \mathcal{H}, D)$ denote a twisted spectral triple obtained by scaling automorphisms of 
an ordinary finitely summable regular spectral triple, introduced in Subsection 
\ref{Twisting spectral triples by scaling automorphisms} and let 
$\Psi (\mathcal{A}_G, \mathcal{H}, D) := \Psi (\mathcal{A}, \mathcal{H}, D) \rtimes G$ denote the associated 
algebra of twisted pseudodifferential operators. The automorphism $\sigma$ extends to an automorphism of 
$\Psi (\mathcal{A}_G, \mathcal{H}, D)$ and it is required that for any $P$ in this algebra \cite{mos1}:
\begin{equation} \label{Selbergcondition}
\int\!\!\!\!\!\!- P = \int\!\!\!\!\!\!- \, \sigma(P).
\end{equation}

\newtheorem{twistedJLO}[first]{Definition}
\begin{twistedJLO} \textnormal{\cite{mos1}}
The twisted JLO bracket of order $q$ as a $q+1$-linear
form on $\Psi (\mathcal{A}_G , \mathcal{H} , D) $ is defined by 
\begin{eqnarray} 
&&   \langle a_0 U_0^* , \dots, a_q U_q^* \rangle_{D} = \int_{\Delta_q} \textnormal{Trace}
  \big(\gamma \, a_0U_0^* \, e^{- s_0 \mu (U_0)^2 D^2} \, a_1 U_1^* \,
  e^{- s_1 \mu (U_0 U_1)^2 D^2} \cdots \nonumber \\
&& \qquad \qquad \qquad   \qquad \qquad   \qquad \qquad
\cdots a_qU_q^* \, e^{- s_q \mu (U_0 \cdots U_q)^2 D^2} \big) \nonumber
\end{eqnarray}
for all $a_0, \dots, a_q \in \mathcal{A}$ and $U_0, \dots, U_q \in G$,
where 
\[
\Delta_q:= \{ s=(s_0, \dots ,s_q)\in \mathbb{R}^{q+1}; \quad s_j\geq 0, \quad s_0+ \cdots +s_q=1 \} . \nonumber
\]
\end{twistedJLO}

In \cite{mos1}, it is assumed that $a_0, \dots, a_q$ are polynomials in $D$ and elements of $\mathcal{A}$
and $[D, \mathcal{A}]$ which are homogeneous in $\lambda$ as $D$ is replaced by $\lambda D$. 
The expression obtained from replacing each
$D$ occurring in $a_0, \dots, a_q$ by $\epsilon^{1/2}D$ is denoted by $\langle a_0 U_0^* , 
\dots, a_q U_q^* \rangle_{D}(\epsilon)$. Equivalently, one has

\begin{eqnarray}
  \langle a_0 U_0^* , \dots, a_q U_q^* \rangle_{D}(\epsilon) \,=\,
  \epsilon^{\frac{m}{2}}   \langle a_0 U_0^* , \dots, a_q U_q^* \rangle_{\epsilon^{1/2} D} \, , \nonumber
\end{eqnarray}
where $m$ is the total degree of $\lambda$ in $a_0 \cdots a_q$ after replacing each $D$ by $\lambda D$.

\newtheorem{asymptotic}[first]{Proposition}
\begin{asymptotic} \textnormal{\cite{mos1}}
If $a_0 \in \mathcal{A}$ and $a_1, \dots, a_q \in [D, \mathcal{A}]$, then there is an asymptotic expansion
\begin{eqnarray}
\langle a_0 U_0^* , \dots, a_q U_q^* \rangle_{D}(\epsilon) \sim \,
   \sum_{j \in J} (c_j + c'_j \log \epsilon)\, \epsilon^{\frac{q}{2}  - \rho_j}  \, + \, O(1) \,\,\,\,\, as \,\,\, \epsilon 
   \searrow 0 \nonumber
\end{eqnarray}
where $\rho_0, \dots, \rho_m$ are points in the half plane $Re(z) \geq \frac{q}{2}. $
\end{asymptotic}

The twisted version of the JLO cocycles is defined by
\begin{equation}
  J^q (D) (A_0,\dots,A_q) \, = \,
  \langle A_0, [D, \sigma^{-1} (A_1)]_{\sigma} , \dots,  [D, \sigma^{-q}(A_q)]_{\sigma} \rangle_{D} \nonumber
\end{equation}
for $A_0, \dots, A_q \in \mathcal{A}_G$. Since in the twisted case this does not define a cocycle, in \cite{mos1}, 
by considering 
\begin{equation}
  J^q (\epsilon^{1/2}D) (A_0,\dots,A_q) \, = \, \epsilon^{\frac{q}{2}}
  \langle A_0, [D, \sigma^{-1} (A_1)]_{\sigma} , \dots,  [D, \sigma^{-q}(A_q)]_{\sigma} 
  \rangle_{\epsilon^{1/2}D} \nonumber
\end{equation}
for $\epsilon > 0$, and passing to the constant term using the above proposition, 
$\varphi_q^\sigma$ is defined by 
\begin{equation}
\varphi_q^\sigma(A_0, \dots, A_q) := \langle A_0, [D, \sigma^{-1} (A_1)]_{\sigma} , 
\dots,  [D, \sigma^{-q}(A_q)]_{\sigma} \rangle_{D} \mid_0. \nonumber
\end{equation}

It follows from the following two results \cite{mos1} that for any twisted spectral triple obtained by 
conformal perturbation of an ordinary finitely summable regular 
spectral triple and satisfying a Selberg type invariance condition \eqref{Selbergcondition}, 
one can associate a local cyclic cocycle for its Connes-Chern character.

\newtheorem{bBcocycle}[first]{Theorem}  
\begin{bBcocycle} \textnormal{\cite{mos1}}
The cochain $\{ \varphi_{q}^\sigma \}$ satisfies the cocycle identity in the (b, B)-bicomplex:
\[b\varphi_{q-1}^\sigma(a_0 U_0^*, \dots, a_q U_q^*) + B\varphi_{q+1}^\sigma(a_0 U_0^*, 
\dots, a_q U_q^*)=0, \nonumber \]
for all $a_0, \dots, a_q \in \mathcal{A}$ and $U_0, \dots, U_q \in G$.
\end{bBcocycle}

\newtheorem{cohom}[first]{Theorem}
\begin{cohom} \textnormal{\cite{mos1}}
The cocycle $\{ \varphi_{q}^\sigma \}$ is cohomologous to the Connes-Chern character 
associated to 
the twisted spectral triple $(\mathcal{A}_G, \mathcal{H}, D)$ in the periodic cyclic 
cohomology $HP^*(\mathcal{A}_G)$.
\end{cohom}

\end{document}